\documentclass[a4paper,10pt]{amsart}

\usepackage{t1enc}

\usepackage[all]{xy}
\usepackage{amssymb}
\usepackage{amsthm}
\usepackage{epic}
\usepackage{amsmath}
\title[On the structure of triangulated categories]{On the structure of triangulated categories with finitely many indecomposables}
\author{Claire Amiot}

\begin{document}

\newcommand{\ten}{\otimes}
\newcommand{\lten}{\overset{L}{\ten}}
\newcommand{\T}{\mathcal{T}}
\newcommand{\Hom}{{\sf Hom }}
\newcommand{\End}{{\sf End }}
\newcommand{\Ext}{{\sf Ext }}
\newcommand{\kd}{k(\mathbb{Z}\Delta)}
\newcommand{\rad}{{\sf rad }}
\newcommand{\modd}{{\sf mod \hspace{.02in} }}
\newcommand{\Aut}{\mbox{Aut}}
\newcommand{\ind}{\mbox{ind}}

\newcommand{\A}{\mathcal{A}cp(\mathcal{E})}

\newcommand{\comp}{{\sf comp }}

\newcommand{\per}{{\sf per \hspace{.02in}}}
\newcommand{\proj}{{\sf proj \hspace{.02in}}}

%
%
\newcommand{\ca}{{\mathcal A}}
\newcommand{\cb}{{\mathcal B}}
\newcommand{\cc}{{\mathcal C}}
\newcommand{\cd}{{\mathcal D}}
\newcommand{\ce}{{\mathcal E}}
\newcommand{\cF}{{\mathcal F}}
\newcommand{\cg}{{\mathcal G}}
\newcommand{\ch}{{\mathcal H}}
\newcommand{\ci}{{\mathcal I}}
\newcommand{\ck}{{\mathcal K}}
\newcommand{\cl}{{\mathcal L}}
\newcommand{\cm}{{\mathcal M}}
\newcommand{\cn}{{\mathcal N}}
\newcommand{\co}{{\mathcal O}}
\newcommand{\cp}{{\mathcal P}}
\newcommand{\cR}{{\mathcal R}}
\newcommand{\cq}{{\mathcal Q}}
\newcommand{\cs}{{\mathcal S}}
\newcommand{\ct}{{\mathcal T}}
\newcommand{\cu}{{\mathcal U}}
\newcommand{\cv}{{\mathcal V}}
\newcommand{\cw}{{\mathcal W}}
\newcommand{\cx}{{\mathcal X}}

\newtheorem{df}{Definition}[subsection]
\newtheorem{lem}{Lemma}[subsection]
\newtheorem{thm}[lem]{Theorem}
\newtheorem{cor}[lem]{Corollary}
\newtheorem{prop}[lem]{Proposition}
\newtheorem{rk}{Remark}

\begin{abstract}
We study the problem of classifying triangulated categories with
finite-dimensional morphism spaces and finitely many indecomposables
over an algebraically closed field. We obtain a new proof of the
following result due to Xiao and Zhu: the Auslander-Reiten quiver of
such a category is of the form $\mathbb{Z}\Delta/G$ where $\Delta$ is
a disjoint union of simply laced Dynkin diagrams and $G$ a weakly
admissible group of automorphisms of $\mathbb{Z}\Delta$. Then we prove
that for `most' groups $G$, the category $\T$ is standard, \emph{i.e.}
$k$-linearly equivalent to an orbit category $\mathcal{D}^b(\modd
k\Delta)/\Phi$. This happens in particular when $\T$ is maximal
$d$-Calabi-Yau with $d\geq2$. Moreover, if $\T$ is standard and
algebraic, we can even construct a triangle equivalence between $\T$ and the corresponding orbit category. Finally we give a sufficient condition for the category of projectives of a Frobenius category to be triangulated. This allows us to construct non standard $1$-Calabi-Yau categories using deformed preprojective algebras of generalized Dynkin type.
\end{abstract}

\maketitle

\section*{Introduction}

Let $k$ be an algebraically closed field and $\T$ a small $k$-linear
triangulated category (see \cite{Ver2}) with split idempotents. We assume that

\hspace{.2in}a) $\T$ is $\Hom$-\emph{finite}, \emph{i.e.} the space
$\Hom_\T(X,Y)$ is finite-dimensional for all objects $X$, $Y$ of $\T$.

It follows that indecomposable objects of $\T$ have local endomorphism
rings and that each object of $\T$ decomposes into a finite direct sum
of indecomposables \cite[3.3]{Gab}. We assume moreover that 

\hspace{.2in}b) $\T$ is \emph{locally finite}, \emph{i.e.} for each
indecomposable $X$ of $\T$, there are at most finitely many isoclasses
of indecomposables $Y$ such that $\Hom_\T(X,Y)\neq 0$.

It was shown in \cite{Xia} that condition $b)$ implies its
dual. Condition $b)$ holds in particular if we have

\hspace{.2in} b') $\T$ is \emph{additively finite}, \emph{i.e.} there
are only finitely many isomorphism classes of indecomposables in $\T$.

The study of particular classes of such triangulated categories $\T$
has a long history. Let us briefly recall some of its highlights:

\hspace{.2in}(1) If $A$ is a representation-finite selfinjective
algebra, then the stable category $\T$ of finite-dimensional (right)
$A$-modules satisfies our assumptions and is additively finite. The
structure of the underlying $k$-linear category of $\T$ was determined
by C. Riedtmann in \cite{Rie}, \cite{Rie2}, \cite{Rie4} and \cite{Rie3}. 
 
\hspace{.2in}(2) In \cite{Hap2}, D. Happel showed that the bounded derived category of
the category of finite-dimensional representations of a representation-finite quiver
is locally finite and described its underlying $k$-linear category.

\hspace{.2in}(3) The stable category $\underline{CM}(R)$ of Cohen-Macaulay modules over a commutative complete local Gorenstein isolated
singularity $R$ of dimension $d$ is a $\Hom$-finite triangulated
category which is $(d-1)$-Calabi-Yau (cf. for
example \cite{Iya} and \cite{Yos}). In \cite{Aus3}, M. Auslander and
I. Reiten showed that if the dimension of R is 1, then the category
$\underline{CM}(R)$ is additively finite and computed the shape of components of its Auslander-Reiten quiver.

\hspace{.2in}(4) The cluster category $\mathcal{C}_Q$ of a finite quiver $Q$ without
oriented cycles was introduced in \cite{Cal} if $Q$ is an orientation
of a Dynkin diagram of type $\mathbb{A}$ and in \cite{Bua} in the
general case. The category $\mathcal{C}_Q$ is triangulated \cite{Kel}
and, if $Q$ is representation-finite, satisfies $a)$ and $b')$.

In a recent article \cite{Xia}, J. Xiao and B. Zhu determined the
structure of the
Auslander-Reiten quiver of a locally finite triangulated category. In
this paper, we obtain the same result with a new proof in section \ref{struc}, namely
that each connected component of the Auslander-Reiten quiver of the
category $\T$ is of the form $\mathbb{Z}\Delta/G$, where $\Delta$ is a
simply laced Dynkin diagram and $G$ is trivial or a weakly admissible
group of automorphisms.
We are interested in the $k$-linear structure of $\T$. If the
Auslander-Reiten quiver of $\T$ is of the form $\mathbb{Z}\Delta$, we
show that the
category $\T$ is standard, \emph{i.e.} it is equivalent to the mesh category $\kd$. Then in section \ref{partcase}, we prove that $\T$ is standard
if the number of vertices of $\Gamma=\mathbb{Z}\Delta/G$ is
strictly greater than the number of isoclasses of
indecomposables of $\modd k\Delta$.
In the last section, using \cite{Bia} we construct examples of non standard triangulated
categories such that $\Gamma = \mathbb{Z}\Delta/\tau$. 

Finally, in the standard cases, we are interested in the triangulated
structure of $\T$. For this, we need to make additional assumptions
on $ \T$. If the Auslander-Reiten quiver is of the form
$\mathbb{Z}\Delta$, and if $\T$ is the base of a tower of triangulated
categories \cite{Kel3}, we show that there is a triangle equivalence between $\T$ and
the derived category $\mathcal{D}^b(\modd k\Delta)$. For the additively finite cases,
we have to assume that $\T$ is standard and algebraic in the sense of \cite{Kel2}. We then show that $\T$ is (algebraically) triangle equivalent to the orbit category of $\mathcal{D}^b(\modd k\Delta)$ under the action of a weakly admissible group of automorphisms. In particular, for each
$d\geq 2$, the algebraic triangulated categories with finitely many
indecomposables which are maximal Calabi-Yau of CY-dimension $d$ are
parametrized by the simply laced Dynkin diagrams.
   
Our results apply in particular to many stable categories
$\underline{\modd}A$ of representation-finite selfinjective algebras
$A$. These algebras were classified up to stable equivalence by
C. Riedtmann \cite{Rie2} \cite{Rie3} and H. Asashiba \cite{Asa}.  In \cite{Bia2}, J. Bia{\l}kowski and A. Skowro\'nski give a
necessary and sufficient condition on these algebras so that
their stable categories $\underline{\modd} A$ are Calabi-Yau. In \cite{Hol}
and \cite{Hol2}, T. Holm and P. J{\o}rgensen prove that certain stable
categories $\underline{\modd}A$ are in fact $d$-cluster categories. These results can also
be proved using our corollary \ref{corcy}. 

This paper is organized as follows: In section \ref{Serre}, we prove that $\T$
has Auslander-Reiten triangles. Section \ref{transquiver} is dedicated to definitions about stable valued translation quivers and admissible
automorphisms groups \cite{Hap3}, \cite{Hap4}, \cite{Die}. We show in
section ~\ref{ART} that the
Auslander-Reiten quiver of $\T$ is a stable valued quiver and in
section ~\ref{struc}, we reprove the result of J. Xiao and B. Zhu \cite{Xia}: The
Auslander-Reiten quiver is a disjoint union of quivers $\mathbb{Z}\Delta/G$,
where $\Delta$ is a Dynkin quiver of type $\mathbb{A}, \mathbb{D}$ or
$\mathbb{E}$, and $G$ a weakly admissible group of automorphisms. In section ~\ref{coveringsection}, we construct a covering functor
$\mathcal{D}^b(\modd k\Delta)\rightarrow \T$ using Riedtmann's method \cite{Rie}.
Then, in section ~\ref{partcase}, we exhibit some combinatorial cases in which $\T$ has to be
standard, in particular when $\T$ is maximal $d$-Calabi-Yau with $d\geq 2$.
Section ~\ref{algsection} is dedicated to the algebraic case. If $\T$ is
algebraic and standard, we can construct a triangle equivalence
between $\T$ and an orbit category.
If $\cp$ is a $k$-category such that $\modd \cp$ is a Frobenius category
satisfying certain conditions, we will prove in section ~\ref{proj} that $\cp$
has naturally a triangulated structure. This allows us to deduce
in section ~\ref{defpreprojalg} that the category $\proj P^f(\Delta)$ of the projective modules over a
deformed preprojective algebra of generalized Dynkin type \cite{Bia}
is naturally triangulated and to reduce the classification of the
additively finite triangulated categories which are $1$-Calabi-Yau to
that of the deformed preprojective algebras in the sense of \cite{Bia}.
 
\section*{Acknowledgments}
I would like to thank my supervisor B. Keller for his availability and many
helpful discussions, and P. J{\o}rgensen for his interest in this work
and for suggesting  some clarifications to me in the proof of theorem \ref{algcy}.  

\section*{Notation and terminology}

We work over an algebraically closed field $k$. By a
\emph{triangulated category}, we mean a $k$-linear triangulated
category $\T$. We
write $S$ for the suspension functor of $\T$ and \begin{math}
\xymatrix{U \ar[r]^u & V \ar[r]^v & W\ar[r]^w & SU}
\end{math} for a distinguished triangle. 
We say that $\T$ is $\Hom$-\emph{finite} if  for each pair $X$, $Y$ of objects in $\T$, the space $\Hom_{\T}(X,Y)$ is 
finite-dimensional over $k$. The category $\T$ will be called a 
\emph{Krull-Schmidt} category if each object is isomorphic to a
finite direct sum of indecomposable objects and the endomorphism ring of an
indecomposable object is a local ring. This implies that idempotents
of $\T$ split \cite[I 3.2]{Hap}.
The category $\T$ will be called
\emph{locally finite} if for each indecomposable $X$ of $\T$,
there are only finitely many isoclasses of indecomposables $Y$
such that $\Hom_\T(X,Y)\neq 0$. This property is selfdual by \cite[prop 1.1]{Xia}.
 
The \emph{Serre functor} will be denoted by $\nu$ (see definition in
section \ref{Serre}). The \emph{Auslander-Reiten translation}
will always be denoted by $\tau$ (section~\ref{Serre}). 

Let $\T$ and $\T'$ be two triangulated categories. An
\emph{$S$-functor} $(F,\phi)$ is given by a $k$-linear functor
$F: \T\rightarrow \T'$ and a functor isomorphism
$\phi$ between the
functors $F\circ S$ and $S'\circ F$, where $S$ is the suspension of
$\T$ and $S'$ the suspension of $\T'$. The notion of $\nu$-functor, or
$\tau$-functor is then clear. A \emph{triangle functor} is an
$S$-functor $(F,\phi)$ such that for each triangle \begin{math}
\xymatrix{U \ar[r]^u & V \ar[r]^v & W\ar[r]^w & SU}
\end{math} of $\T$, the sequence \begin{math}
\xymatrix{FU \ar[r]^{Fu} & FV \ar[r]^{Fv} & FW\ar[rr]^{\phi_{U}\circ
    Fw} & &S'FU}
\end{math} is a triangle of $\T'$.

The category $\T$ is \emph{Calabi-Yau} if there exists an integer $d>0$ such
that we have a triangle functor isomorphism between $S^d$ and $\nu$. We say
that $\T$ is \emph{maximal d-Calabi-Yau} if $\T$ is $d$-Calabi-Yau and
if for each covering functor
$\T'\rightarrow \T$ with $\T'$ $d$-Calabi-Yau, we have a $k$-linear
equivalence between $\T$ and $\T'$.

For an additive $k$-category $\mathcal{E}$, we write $\modd\mathcal{E}$ for
the category of contravariant finitely presented functors from
$\mathcal{E}$ to $\modd k$ (section \ref{proj}), and if the
projectives of $\modd \mathcal{E}$ coincide with the
injectives, $\underline{\modd}\mathcal{E}$ will be the \emph{stable category}. 

\section{Serre duality and Auslander-Reiten triangles}\label{Serre}

\subsection{Serre duality}

Recall from \cite{Rei} that a \emph{Serre functor} for $\T$ is an
autoequivalence $\nu:\T\rightarrow \T$ together with an isomorphism
$D\Hom_\T(X,?)\simeq \Hom_\T(?,\nu X)$ for each $X\in\T$, where $D$
is the duality $\Hom_k(?,k)$.
 
\begin{thm}
Let $\T$ be a Krull-Schmidt, locally finite triangulated category. Then $\T$
has a Serre functor $\nu$.
\end{thm}
\begin{proof}
Let $X$ be an object of $\T$. We write $X^\wedge$ for the functor
$\Hom_\T(?,X)$ and $F$ for the functor $D\Hom_{\T}(X,?)$. Using the lemma \cite[I.1.6]{Rei} we just have to
show that $F$ is representable. Indeed, the category $\T^{op}$ is
locally finite as well.
The proof is in two steps.

\emph{Step 1: The functor
  $F$ is finitely presented.}

Let $Y_1,\ldots, Y_r$ be representatives of the isoclasses of indecomposable
objects of $\T$ such that $FY_i$ is not zero. The space $\Hom(Y_i^\wedge,F)$ is
finite-dimensional over $k$. Indeed it is isomorphic to $FY_i$ by the
Yoneda lemma. Therefore, the functor
 $\Hom(Y_i^\wedge,F)\ten_k Y_i^\wedge$ is representable.
We get an epimorphism from a representable functor to $F$:
$$\bigoplus_{i=1}^r \Hom(Y_i^\wedge,F)\ten_k Y_i^\wedge \longrightarrow F.$$
By applying the same argument to its kernel we get a
projective presentation of $F$ of the form
$U^\wedge \longrightarrow V^\wedge \longrightarrow F \longrightarrow
0$, with $U$ and $V$ in $\T$. 

\emph{Step 2: A cohomological functor}
$H:\T^{op}\rightarrow \modd  k$ \emph{is representable  if and only if it is
finitely presented}.

Let 
\begin{math}
\xymatrix{U^\wedge \ar[r]^{u^\wedge} & V^\wedge \ar[r]^\phi & H \ar[r]  & 0}
\end{math} be a presentation of $H$. We form a triangle  
\begin{math}
\xymatrix{U \ar[r]^u & V \ar[r]^v & W\ar[r]^w & SU.}
\end{math}
We get an exact sequence 
$$\xymatrix{U^\wedge \ar[r]^{u^\wedge} & V^\wedge \ar[r]^{v^\wedge} & W^\wedge \ar[r]^{w^\wedge} & (SU)^{\wedge}}
.$$ Since the composition of $\phi$ with $u^\wedge$ is zero and $H$ is
cohomological, the morphism
$\phi$ factors through $v^\wedge$. But $H$ is the cokernel of
$u^\wedge$, so $v^\wedge$ factors through $\phi$.
We obtain a commutative diagram:
\begin{displaymath}
\xymatrix{U^\wedge \ar[r]^{u^\wedge} & V^\wedge \ar[r]^{v^\wedge}\ar[d]_\phi & W^\wedge \ar[r]^{w^\wedge}\ar@<1ex>[dl]^{\phi'} & SU.\\
            & H\ar@<1ex>[ur]^i & &}
\end{displaymath}
The equality $\phi'\circ i \circ
\phi=\phi'\circ v^\wedge=\phi$ implies that $\phi'\circ i$ is the
identity of $H$ because $\phi$ is an epimorphism. We deduce that $H$ is
a direct factor of $W^\wedge$. The composition $i\circ \phi'=e^\wedge$
is an idempotent. Then $e\in \End(W)$ splits and we get  $H=W'^\wedge$
for a direct factor $W'$ of $W$.

\end{proof} 

\subsection{Auslander-Reiten triangles}

\begin{df}
\cite{Hap2} A triangle  
$\xymatrix{ X \ar^u[r] & Y\ar^v[r] & Z\ar^w[r] & SX}$
of $\ \T$ is called an \emph{Auslander-Reiten triangle} or \emph{AR-triangle} if the following conditions
 are satisfied:

\hspace{.2in} (AR1) $X$ and $Z$ are indecomposable objects;

\hspace{.2in} (AR2) $w\neq 0$;

\hspace{.2in} (AR3) if $f:W\rightarrow Z$ is not a retraction, there exists
    $f':W\rightarrow Y$ such that $vf'=f$;
 
\hspace{.2in} (AR3') if $g:X\rightarrow V$ is not a section, there exists
$g':Y\rightarrow V$ such that $g'u=g$.
\end{df}
Let us recall that, if (AR1) and (AR2) hold, the conditions (AR3) and
(AR3') are equivalent.
We say that a triangulated category $\T$ \emph{has Auslander-Reiten
triangles} if, for any indecomposable object $Z$ of $\T$, there exists an
AR-triangle ending at $Z$:
$\xymatrix{ X \ar^u[r] & Y\ar^v[r] & Z\ar^w[r] & SX}$.
In this case, the AR-triangle is unique up to triangle isomorphism
inducing the identity of $Z$.

The following proposition is proved in \cite[Proposition I.2.3]{Rei}
\begin{prop}
The category $\T$ has Auslander-Reiten triangles.
\end{prop}
The composition $\tau= S^{-1}\nu$ is called
the Auslander-Reiten translation. An AR-triangle of $\T$ ending at $Z$
has the form:
$$\xymatrix{ \tau Z \ar^u[r] & Y\ar^v[r] & Z\ar^w[r] & \nu Z.}$$

\section{Valued translation quivers and automorphism groups}\label{transquiver}

\subsection{Translation quivers}
In this section, we recall some definitions and notations concerning
quivers \cite{Die}.
A quiver $Q=(Q_0, Q_1, s, t)$ is given by the set $Q_0$ of its
vertices, the set $Q_1$ of its arrows, a source map $s$ and
a tail map $t$. If $x\in Q_0$ is a vertex, we denote by $x^+$ the set
of direct successors of $x$, and by $x^-$ the set of its direct
predecessors. We say that $Q$ is \emph{locally finite} if for each
vertex $x\in Q_0$, there are finitely many arrows ending at $x$ and
starting at $x$ (in this case, $x^+$ and $x^-$ are finite sets). The quiver $Q$ is said to be \emph{without double arrows},
if two different arrows cannot have the same tail and source.

\begin{df}
A \emph{stable translation quiver} $(Q,\tau)$ is a locally finite quiver without
double arrows with a bijection $\tau:Q_0\rightarrow Q_0$ such that
$(\tau x)^+=x^-$ for each vertex $x$.
For each arrow $\alpha :x\rightarrow y$, let $\sigma\alpha$ the unique arrow
$\tau y \rightarrow x$. \end{df} 
Note that a stable translation quiver can have loops.

\begin{df} 
A \emph{valued translation quiver} $(Q,\tau ,a)$ is a stable
translation quiver $(Q,\tau )$ with a map $a:Q_1\rightarrow
\mathbb{N}$ such that $a(\alpha)=a(\sigma\alpha)$ for each arrow
$\alpha$.
If $\alpha$ is an arrow from $x$ to $y$, we write $a_{xy}$ instead of $a(\alpha)$.
\end{df} 
\begin{df}
Let $\Delta$ be an oriented tree. The \emph{repetition of } $\Delta$
is the quiver $\mathbb{Z}\Delta$ defined as follows:
\begin{itemize}
\item $(\mathbb{Z}\Delta)_0=\mathbb{Z}\times\Delta_0$
\item $(\mathbb{Z}\Delta)_1=\mathbb{Z}\times \Delta_1 \cup \sigma(
  \mathbb{Z}\times \Delta_1)$ with arrows
  $(n,\alpha):(n,x)\rightarrow (n,y)$ and
  $\sigma(n,\alpha):(n-1,y)\rightarrow (n,x)$ for each arrow $\alpha
  :x\rightarrow y$ of $\Delta$.
\end{itemize}
\end{df}
The quiver $\mathbb{Z}\Delta$ with the translation $\tau(n,x)=(n-1,x)$
is clearly a stable translation quiver which does not depend (up to
isomorphism) on the
orientation of $\Delta$ (see \cite{Rie}).

\subsection{Groups of weakly admissible automorphisms}\label{groupe}

\begin{df}
An automorphism group $G$ of a quiver is said to be \emph{admissible \cite{Rie}}
if no orbit of $G$ intersects a set of the form $\{x\}\cup x^+$ or
$\{x\}\cup x^-$ in more than one 
point. It said to be
\emph{weakly admissible \cite{Die}} if, for each $g\in G-\{1\}$ and for
each $x\in Q_0$, we have $x^+\cap (gx)^+=\emptyset$.
\end{df}
Note that an admissible automorphism group is a weakly admissible
automorphism group.
Let us fix a numbering and an orientation of the simply-laced Dynkin trees.

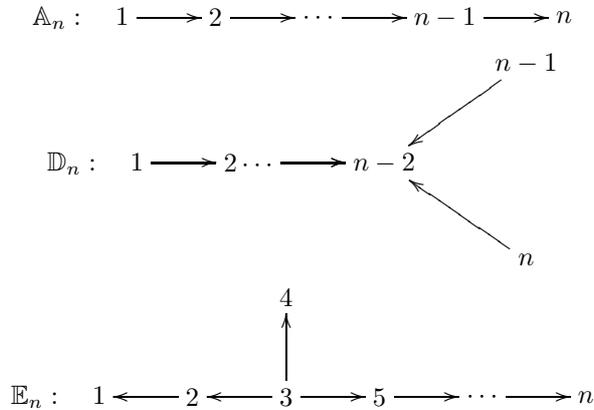
\begin{displaymath}
\xymatrix{\mathbb{A}_n:\quad 1\ar[r] & 2\ar[r] & \cdots \ar[r]& n-1\ar[r] & n}
\end{displaymath}
$$\xymatrix{& & &  n-1\ar[dl]\\ \mathbb{D}_n:\quad 1\ar[r] & 2\ar[r] \cdots\ar[r] & n-2 & \\ & & & n\ar[ul]}$$
$$\xymatrix{ & & 4 & & & \\ \mathbb{E}_n:\quad 1 & 2\ar[l] &
  3\ar[l]\ar[r]\ar[u] & 5\ar[r] & \cdots \ar[r] & n}$$

Let $\Delta$ be a Dynkin tree. We define an automorphism $S$ of
$\mathbb{Z}\Delta$ as follows:
\begin{itemize}
\item if $\Delta=\mathbb{A}_n$, then $S(p,q)=(p+q,n+1-q)$;
\item if $\Delta=\mathbb{D}_n$ with $n$ even, then $S=\tau^{-n+1}$;
\item if $\Delta=\mathbb{D}_n$ with $n$ odd, then
  $S=\tau^{-n+1}\phi$ where $\phi$ is the automorphism of $\mathbb{D}_n$
 which exchanges $n$ and $n-1$;
\item if $\Delta=\mathbb{E}_6$, then $S=\phi\tau^{-6}$ where $\phi$ is
  the automorphism of $\mathbb{E}_6$ which exchanges $2$ and $5$,
  and $1$ and $6$;
\item if $\Delta=\mathbb{E}_7$, then $S=\tau^{-9}$;
\item and if $\Delta=\mathbb{E}_8$, then $S=\tau^{-15}$.
\end{itemize}

In \cite[Anhang 2]{Rie}, Riedtmann describes all admissible
automorphism groups of Dynkin diagrams. Here is a more precise result:
\begin{thm}\label{groupesadmissibles}
Let $\Delta$ be a Dynkin tree and $G$ a non trivial group of weakly admissible
automorphisms of $\mathbb{Z}\Delta$. Then $G$ is isomorphic to
$\mathbb{Z}$, and here is a list of its possible generators:

\begin{itemize}
\item
if $\Delta=\mathbb{A}_n$ with $n$ odd, possible generators are
$\tau^r$ and $\phi\tau^r$ with $r\geq 1$, where $\phi=\tau^\frac{n+1}{2}S$
is an automorphism of $\mathbb{Z}\Delta$ of order 2;
\item
if $\Delta=\mathbb{A}_n$ with $n$ even, then possible generators are $\rho^r$, where $r\geq 1$ and
where $\rho=\tau^{\frac{n}{2}}S$. (Since $\rho^2=\tau^{-1}$, $\tau^r$
is a possible generator.)
\item
if $\Delta=\mathbb{D}_n$ with $n\geq 5$, then possible generators are
$\tau^r$ and $\tau^r\phi$, where $r\geq 1$ and
where $\phi= (n-1,n)$ is the automorphism of $\mathbb{D}_n$ exchanging $n$
and $n-1$.
\item
if $\Delta=\mathbb{D}_4$, then possible generators are
$\phi\tau^r$, where $r\geq 1$ and where $\phi$ belongs to
$\mathfrak{S}_3$ the permutation group on 3 elements seen as subgroup of
automorphisms of $\mathbb{D}_4$.
\item
if $\Delta=\mathbb{E}_6$, then possible generators are
$\tau^r$ and $\phi\tau^r$, where $r\geq 1$ and where $\phi$ is the
automorphism of $\mathbb{E}_6$ exchanging $2$ and $5$, and $1$ and $6$.
\item
if $\Delta=\mathbb{E}_n$ with $n=7,8$, possible generators are $\tau^r$, where $r\geq 1$.   
\end{itemize}
The unique weakly admissible automorphism group which is not
admissible exists for $\mathbb{A}_n$, $n$ even, and is generated by $\rho$.
\end{thm}

\section{Property of the Auslander-Reiten translation }\label{ART}

We define the Auslander-Reiten quiver $\Gamma_\T$ of the category $\T$
as a valued quiver $(\Gamma,a)$. The vertices are the isoclasses of
indecomposable objects. Given two indecomposable objects $X$ and $Y$
of $\T$, we draw one arrow from $x=[X]$ to $y=[Y]$ if the vector
space $\mathcal{R}(X,Y)/\mathcal{R}^2(X,Y)$ is not zero, where
$\mathcal{R}(?,?)$ is the radical of the bifunctor $\Hom_\T(?,?)$.
A morphism of $\mathcal{R}(X,Y)$ which does not vanish in the quotient
$\mathcal{R}(X,Y)/\mathcal{R}^2(X,Y)$ will be called
\emph{irreducible}. Then we put 
 $$a_{xy}=\dim_k\mathcal{R}(X,Y)/\mathcal{R}^2(X,Y).$$

Remark that the fact that $\T$ is locally finite implies that its AR-quiver is locally finite. The aim of this section is to show that $\Gamma_\T$ with the
translation $\tau$ defined in the first part is a valued translation quiver. In other words, we want to show the proposition:

\begin{prop}\label{trans}
If $X$ and $Y$ are indecomposable objects of $\T$, we have the equality 
$$\dim_k \mathcal{R} (X,Y)/ \mathcal{R}^2(X,Y)=\dim_k \mathcal{R} (\tau Y,X)/ \mathcal{R}^2(\tau Y,X).$$
\end{prop}

Let us recall some definitions \cite{Hap}.
\begin{df}
A morphism $g:Y\rightarrow Z$ is called \emph{sink morphism} if
the following hold

(1) $g$ is not a retraction;

(2) if $h:M\rightarrow Z$ is not a retraction, then $h$ factors
    through $g$;

(3) if $u$ is an endomorphism of $Y$ which satisfies $gu=u$, then $u$
    is an automorphism.

Dually, a morphism $f:X\rightarrow Y$ is called  \emph{source
  morphism} if the following hold:

(1) $f$ is not a section;

(2) if $h:X\rightarrow M$ is not a section, then $h$ factors through
    $f$;

(3) if $u$ is an endomorphism of $Y$ which satisfies $uf=f$, then $u$
    is an automorphism.

\end{df}
These conditions imply that $X$ and $Z$ are indecomposable.
Obviously, if \\ $\xymatrix{X\ar[r]^u & Y\ar[r]^v & Z\ar[r]^w
  & SX}$ is an AR-triangle, then $u$ is a source morphism and $v$ is a
  sink morphism. Conversely, if  $v\in \Hom_\T(Y,Z)$ is a sink
  morphism (or if $u\in\Hom_\T(X,Y)$ is a source morphism), then there
  exists an AR-triangle $\xymatrix{X\ar[r]^u & Y\ar[r]^v & Z\ar[r]^w &
  SX}$ (see \cite[I 4.5]{Hap}).

The following lemma (and the dual statement) is proved in \cite[2.2.5]{Rin}.

\begin{lem}
Let $g$ be a morphism from $Y$ to $Z$, where $Z$ is indecomposable
and 
 $Y=\bigoplus_{i=1}^r Y_i^{n_i}$ is the decomposition of $Y$ into
 indecomposables. Then the morphism $g$ is a sink morphism if and only
 if the following hold:

(1) For each $i=1,\ldots, r$ and $j=1,\ldots,n_i$, the morphism
    $g_{i,j}$ belongs to the radical $\mathcal{R}(Y_i,Z)$.

(2) For each $i=1,\ldots, r$, the family
    $(\overline{g}_{i,j})_{j=1,\ldots,n_i}$ forms a $k$-basis of the space $\mathcal{R}(Y_i,Z)/\mathcal{R}^2(Y_i,Z)$.

(3) If $h\in \Hom_\T(Y',Z)$ is irreducible and $Y'$ indecomposable,
    then $h$ factors through $g$ and $Y'$ is isomorphic to $Y_i$
    for some $i$.

\end{lem}

Using this lemma, it is easy to see that proposition \ref{trans}
holds. Thus, the Auslander-Reiten quiver $\Gamma_\T =(\Gamma,
\tau, a)$ of the
category $\T$ is a valued translation quiver.
 
\section{Structure of the Auslander-Reiten quiver}\label{struc}

This section is dedicated to an other proof of a theorem due to J. Xiao
and B. Zhu (\cite{Xia2}):

\begin{thm}\label{structure}\cite{Xia2} 
Let $\T$ be a Krull-Schmidt, locally finite triangulated category.
Let $\Gamma$ be a connected component of the AR-quiver of $\T$. Then
there exists a Dynkin tree $\Delta$ of type $\mathbb{A}$,
$\mathbb{D}$ or $\mathbb{E}$, a weakly admissible automorphism group
$G$ of $\mathbb{Z}\Delta$ and an isomorphism of valued translation quiver 
$$\xymatrix@1{\theta:\Gamma\ar[r]^(.45){\sim} &\mathbb{Z}\Delta/G}.$$
The underlying graph of the tree $\Delta$ is unique up to isomorphism
(it is called the \emph{type of $\Gamma$}), and the group $G$ is unique
up to conjugacy in $\Aut(\mathbb{Z}\Delta)$.

In particular, if $\T$ has an infinite number of isoclasses of
indecomposable objects, then $G$ is trivial, and $\Gamma$ is the
repetition quiver $\mathbb{Z}\Delta$.
\end{thm}
  
\subsection{Auslander-Reiten quivers with a loop}

In this section, we suppose that the Auslander-Reiten quiver of $\T$
contains a loop, \textit{i.e.} there exists an arrow with same tail and source. Thus,
we suppose that there exists an indecomposable $X$ of $\T$ such that
$$\dim_k\mathcal{R}(X,X)/\mathcal{R}^2(X,X)\geq 1.$$ 
\begin{prop}\label{boucle}
Let $X$ be an indecomposable object of $\T$. Suppose that we have  
$\dim_k\mathcal{R}(X,X)/\mathcal{R}^2(X,X)\geq 1.$ Then $\tau X$ is
isomorphic to $X$.
\end{prop}
To prove this, we need a lemma.
\begin{lem}\label{boucle2}
Let \begin{math}\xymatrix{X_1\ar[r]^{f_1} & X_2 \ar[r]^{f_2} &
    \cdots\ar[r]^{f_n} & X_{n+1}}\end{math} be a sequence of irreducible
    morphisms between indecomposable objects with $n\geq 2$. If the composition
    $f_n\circ f_{n-1}\cdots f_1$ is zero, then there exists an $i$ such that
    $\tau ^{-1}X_i$ is isomorphic to $X_{i+2}$.
\end{lem}
\begin{proof}The proof proceeds by induction on $n$. Let us show the assertion
    for $n=2$. Suppose $\xymatrix{X_1\ar[r]^{f_1} &
    X_2\ar[r]^{f_2} & X_3}$ is a sequence such that $f_2\circ f_1=0$. We can then
    construct an AR-triangle:
\begin{displaymath}
\xymatrix{X_1\ar[r]^{(f_1,f)^T} & X_2\oplus X
  \ar[r]^{(g_1,g_2)}\ar[d]^{(f_2,0)}
  &\tau^{-1}X_1\ar[r]\ar@{.>}[dl]^\beta & SX_1\\ & X_3 & & }\end{displaymath}
 The composition $f_2\circ f_1$ is zero, thus the morphism $f_2$ factors
  through $g_1$. As the morphisms $g_1$ and $f_2$ are irreducible, we
  conclude that $\beta$ is a retraction, and $X_3$ a direct summand of
  $\tau^{-1}X_1$. But $X_1$ is indecomposable, so $\beta$ is an
  isomorphism between $X_3$ and $\tau^{-1}X_1$.

Now suppose that the property holds for an integer $n-1$ and that we
have $f_n f_{n-1} \cdots f_1=0$. If the composition $f_{n-1}\cdots
f_1$ is zero, the proposition holds by induction. So we can suppose that
for $i\leq n-2$, the objects $\tau^{-1} X_i$ and $X_{i+2}$ are not isomorphic. We
show now by induction on $i$ that for each $i\leq n-1$, there exists a
map $\beta_i:\tau^{-1}X_i\rightarrow X_{n+1}$ such that $f_n\cdots
f_{i+1}=\beta_ig_i$ where $g_i:X_{i+1}\rightarrow \tau^{-1}X_i$ is an irreducible morphism.
For $i=1$, we construct an AR-triangle: 
\begin{displaymath}
\xymatrix{X_1\ar[r]^{(f_1,f'_1)^T} & X_2\oplus X'_1
                      \ar[r]^{(g_1,g'_1)}\ar[d]_{(f_n\cdots f_2,0)} 
                      &\tau^{-1}X_1\ar[r]\ar@{.>}[dl]^{\beta_1} & SX_1\\
                      & X_{n+1} & & }\end{displaymath}
As the composition $f_n\cdots f_1$ is zero, we have the factorization
                      $f_n\cdots f_2=\beta _1 g_1$.

Now for $i$, as $\tau^{-1} X_{i-1}$ is not isomorphic to $X_{i+1}$, there
exists an AR-triangle of the form:
\begin{displaymath}
\xymatrix@1{X_i\ar^(.3){(g_{i-1}, f_i,f'_i)^T}[rr] & & \tau^{-1}X_{i-1}\oplus
                      X_{i+1}\oplus X'_i \ar^(.6){(g''_i,g_i,g'_i)}[rr]\ar[d]_{(-\beta_{i-1},f_n\cdots
                      f_{i+1},0)} & &\tau^{-1}X_i\ar[r]\ar@{.>}[dll]^{\beta_i} & SX_i\\
                      & & X_{n+1} & & & }\end{displaymath}
By induction,  $-\beta_{i-1}g_{i-1}
                      +f_n\cdots f_{i+1}f_i$ is zero, thus $f_n\cdots
                      f_{i+1}$ factors through $g_i$. This property is
                      true for $i=n-1$, so we have a map
                      $\beta_{n-1}:\tau^{-1}X_{n-1}\rightarrow
                      X_{n+1}$ such that $\beta_{n-1}g_{n-1}=f_n$. As $g_{n-1}$ and
                      $f_n$ are irreducible, we conclude that
                      $\beta_{n-1}$ is an isomorphism between
                      $X_{n+1}$ and $\tau^{-1}X_{n-1}$.
 \end{proof} 

Now we are able to prove proposition \ref{boucle}. There exists an irreducible map $f:X\rightarrow
X$. Suppose that $X$ and $\tau X$ are not isomorphic. Then from the
previous lemma, the endomorphism $f^n$ is non
zero for each $n$. But since $\T$ is a Krull-Schmidt, locally finite
category,  a power of the radical
$\mathcal{R}(X,X)$ vanishes. This is a
contradiction.

\subsection{Proof of theorem \ref{structure}}
Let 
$\tilde{\Gamma}=(\tilde{\Gamma}_0,\tilde{\Gamma}_1, \tilde{a})$ be the
valued translation quiver obtained from $\Gamma$ by removing the loops, \textit{i.e.} we have
$\tilde{\Gamma}_0=\Gamma_0$, $\tilde{\Gamma}_1=\{\alpha\in\Gamma_1$
such that $s(\alpha )\neq t(\alpha )\}$, and $\tilde{a}=a_{|_{\tilde{\Gamma}_1}}$.  
                                            
\begin{lem}
The quiver $\tilde{\Gamma}=(\tilde{\Gamma}_0,\tilde{\Gamma}_1,
\tilde{a})$ with the translation $\tau$ is a valued translation quiver
without loop. 
\end{lem}

\begin{proof}
We have to check that the map $\sigma$ is well-defined. But from
proposition \ref{boucle}, if $\alpha$ is a loop on a vertex $x$,
$\sigma(\alpha)$ is the unique arrow from $\tau x=x$ to $x$,
\emph{i.e.} $\sigma(\alpha)=\alpha$. Thus $\tilde{\Gamma}$ is obtained  from $\Gamma$
by removing some $\sigma$-orbits and it keeps the structure of stable
valued translation quiver.
\end{proof}   

Now, we can apply Riedtmann's Struktursatz \cite{Rie} and the result of Happel-Preiser-Ringel \cite{Hap4}. There exist a tree
$\Delta$ and an admissible automorphism group $G$ (which may be
trivial) of
$\mathbb{Z}\Delta$ such that $\tilde{\Gamma}$ is isomorphic to
$\mathbb{Z}\Delta /G$ as a valued translation quiver. The underlying graph
of the tree $\Delta$ is then unique up to isomorphism and the group
$G$ is unique up to conjugacy in $\Aut(\mathbb{Z}\Delta)$.
Let $x$ be a vertex of $\Delta$. We write $\overline{x}$ for the image of $x$ by the map:       
$$\xymatrix@1{\Delta\ar[r] &\mathbb{Z}\Delta\ar[r]^(.3)\pi &
  \mathbb{Z}\Delta/G\simeq\tilde{\Gamma}\ \ar@{^{(}->}[r] & \Gamma}.$$ 

Let $C:\Delta_0\times\Delta_0\rightarrow \mathbb{Z}$ be the matrix
defined as follows: 
\begin{itemize}
\item
$C(x,y)=-a_{\overline{x}\, \overline{y}}$ (resp. $-a_{\overline{y}\,
  \overline{x}}$) if there exists an arrow from 
$x$ to $y$ (resp. from $y$ to
$x$) in $\Delta$,
\item $C(x,x)=2-a_{\overline{x}\, \overline{x}}$,
\item $C(x,y)=0$ otherwise.
\end{itemize} 
The matrix $C$ is symmetric; it is a `generalized Cartan matrix' in
the sense of \cite{Hap3}. If we remove the loops from the `underlying graph of
$C$' (in the sense of \cite{Hap3}), we get the underlying graph of $\Delta$.

In order to apply the result of Happel-Preiser-Ringel \cite[section
2]{Hap3}, we have to show:

\begin{lem}
The set $\Delta_0$ of vertices of $\Delta$ is finite.
\end{lem}

\begin{proof}
Riedtmann's construction of $\Delta$ is the following. We fix a vertex
$x_0$ in $\tilde{\Gamma}_0$. Then the vertices of $\Delta$ are the
paths of $\tilde{\Gamma}$ beginning on $x_0$ and which do not contain
subpaths of the form $\alpha \sigma(\alpha)$, where $\alpha$ is in
$\tilde{\Gamma}_1$.
Now suppose that $\Delta_0$ is an infinite set. Then for each $n$, there exists a
sequence:
$$\xymatrix{x_0\ar[r]^{\alpha_1} & x_1\ar[r]^{\alpha_2} & \cdots
  \ar[r]^{\alpha_{n-1}} & x_{n-1}\ar[r]^{\alpha_n}& x_n }$$
such that $\tau x_{i+2}\neq x_i$. Then there exist some
  indecomposables $X_0,\ldots, X_n$ such that the vector space
  $\mathcal{R}(X_{i-1}, X_i)/\mathcal{R}^2(X_{i-1}, X_i)$ is not
  zero. Thus from the lemma \ref{boucle2},
  there exists irreducible morphisms $f_i:X_{i-1}\rightarrow X_i$ such
  that the composition $f_n f_{n-1}\cdots f_1$ does not vanish. But
  the functor $\Hom_\T(X_0,?)$ has finite support. Thus there is an
  indecomposable $Y$ which appears an infinite number of times in the
  sequence $(X_i)_i$. But since $\mathcal{R}^N(Y,Y)$ vanishes for an
  $N$, we have a contradiction.  
\end{proof}

Let $\mathcal{S}$ a system of representatives of isoclasses of indecomposables
of $\T$.
For an indecomposable $Y$ of $\T$, we put $$l(Y)=\sum_{M\in
  \mathcal{S}}\dim_k \Hom_\T(M,Y).$$
This sum is finite since $\T$ is locally finite.

\begin{lem}
For $x$ in $\Delta_0$, we write $d_x=l(\overline{x})$. Then for each
$x\in\Delta_0$, we have: $$\sum_{y\in \Delta_0} d_y C_{xy}=2.$$
\end{lem} 
\begin{proof}
Let $X$ and $U$ be indecomposables of $\T$. Let
$$\xymatrix{X\ar[r]^u & Y\ar[r]^v & Z\ar[r]^w & SX}$$ 
be an AR-triangle. We write $(U,?)$
for the cohomological functor $\Hom_\T(U,?)$. Thus, we have a long
exact sequence:
$$\xymatrix@1{(U,S^{-1}Z)\ar[r]^(.6){S^{-1}w_*} & (U,X)\ar[r]^{u_*}
  &(U,Y)\ar[r]^{v_*} & (U,Z)\ar[r]^{w_*} &(U,SX)}.$$
Let  $S_Z(U)$ be the image of the map $w_*$. We have the exact
  sequence: 
$$\xymatrix{0\ar[r] & S_{S^{-1}Z} (U)\ar[r] &  (U,X)\ar[r]^{u_*} &
  (U,Y) \ar[r]^{v_*} 
 &(U,Z)\ar[r]^{w_*} & S_Z(U)\ar[r] & 0.}$$
Thus we have the following equality:
$$\dim_k S_Z(U) + \dim_k S_{S^{-1}Z}(U) + \dim_k (U,Y)=\dim_k(U,X)+\dim_k(U,Z).$$
If $U$ is not isomorphic to $Z$, each map from $U$ to $Z$ is radical,
thus $S_Z(U)$ is zero. If $U$ is isomorphic to $Z$, the map $w_*$
factors through the radical of $\End(Z)$, so $S_Z(Z)$ is isomorphic
to $k$. Then summing the previous equality when $U$ runs over $\mathcal{S}$, we get:
$$l(X)+l(Z)=l(Y)+2.$$ 
Clearly $l$ is $\tau$-invariant, thus $l(Z)$ equals $l(X)$. If the decomposition of $Y$
is of the form $\bigoplus_{i=1}^rY_i^{n_i}$, we get: 
$$l(Y)=\sum_in_i l(Y_i)=\sum_{i, X\rightarrow Y_i \in\tilde{\Gamma}} a_{XY_i} l(Y_i)+ a_{XX} l(X).$$
We deduce the formula:
$$2=(2-a_{XX})l(X) - \sum_{i, X\rightarrow Y_i \in\tilde{\Gamma}} a_{XY_i} l(Y_i).$$
Let $x$ be a vertex of the tree $\Delta$ and $\overline{x}$ its image
in $\tilde{\Gamma}$. Then an arrow $\overline{x}\rightarrow Y$ in
$\tilde{\Gamma}$ comes from an arrow $(x,0)\rightarrow (y,0)$ in
$\mathbb{Z}\Delta$ or from
an arrow $(x,0)\rightarrow (y,-1)$ in $\mathbb{Z}\Delta$, \textit{i.e.} from an
arrow $(y,0)\rightarrow (x,0)$. Indeed the
projection $\mathbb{Z}\Delta\rightarrow \mathbb{Z}\Delta /G$ is a
covering. From this we deduce the following equality:

$$2=(2-a_{\overline{x}\,\overline{x}})d_x -\sum_{y,x\rightarrow y
  \in\Delta} a_{\overline{x}\,\overline{y}} d_y -\sum_{y,y\rightarrow
  x \in\Delta} a_{\overline{y}\,\overline{x}} d_y=\sum_{y\in \Delta_0} d_y C_{xy}.$$ 
\end{proof}

Now we can prove theorem~\ref{structure}. The matrix $C$
is a `generalized Cartan matrix'. The previous lemma gives us a
subadditive function which is not additive. Thus by \cite{Hap3}, the
underlying graph of $C$ is of `generalized Dynkin type'. As $C$ is
symmetric, the graph is necessarily of type $\mathbb{A}$, $\mathbb{D}$,
$\mathbb{E}$, or $\mathbb{L}$. But this graph is the graph $\Delta$ with
the valuation $a$. We are done in the cases $\mathbb{A}$,
$\mathbb{D}$, or $\mathbb{E}$.

The case $\mathbb{L}_n$ occurs when the AR-quiver
contains at least one loop. We can see
$\mathbb{L}_n$ as $\mathbb{A}_n$ with valuations on the vertices with
a loop. Then, it is obvious that the automorphism groups of
$\mathbb{ZL}_n$ are generated by $\tau^r$ for an $r\geq 1$. But
proposition \ref{boucle} tell us that a vertex $x$ with a loop
satisfies $\tau x =x$. Thus $G$ is
generated by $\tau$ and the AR-quiver has the following form:  
$$\xymatrix{1\ar@/^/[r] & 2\ar@/^/[l]\ar@/^/[r] & 3\ar@/^/[l] \ar@{.}[r]& 
  \ar@/^/[r] & n \ar@/^/[l] \ar@(ur,dr)}$$
This quiver is isomorphic to the quiver $\mathbb{ZA}_{2n}/G$ where
  $G$ is the group generated by the automorphism $\tau^{n}S=\rho$.

The suspension functor $S$ sends the indecomposables on
  indecomposables, thus it can be seen as an automorphism of the
  AR-quiver. It is exactly the automorphism $S$ defined in section
  \ref{groupe}. 

As shown in \cite{Xia2}, it follows from the results of \cite{Kel} that for each Dynkin tree $\Delta$ and for each weakly admissible group of
automorphisms $G$ of $\mathbb{Z}\Delta$,  there exists a locally finite triangulated category $\T$ such that 
$\Gamma_\T\simeq \mathbb{Z}\Delta/G$. This category is of the form $\T=\mathcal{D}^b(\modd k\Delta)/\varphi$ where $\varphi$ is an auto-equivalence of $\mathcal{D}^b(\modd k\Delta)$.

\section{Construction of a covering functor}\label{coveringsection}

From now, we suppose that the AR-quiver $\Gamma$ of $\T$ is
connected. We know its structure. It is natural
to ask: Is the category $\T$ \emph{standard}, \emph{i.e.} equivalent
as a $k$-linear category to
the mesh category $k(\Gamma)$? First, in this part we construct a covering functor
$F:k(\mathbb{Z}\Delta)\rightarrow \T$.

\subsection{Construction}
We write $\pi:\mathbb{Z}\Delta\rightarrow \Gamma$ for the
canonical projection. As $G$ is a weakly admissible group, this projection verifies the following property:
if $x$ is a vertex of $\mathbb{Z}\Delta$, the number of
arrows of $\mathbb{Z}\Delta$ with source $x$ is equal to the
number of arrows of $\mathbb{Z}\Delta /G$ with source $\pi x$.
Let $\mathcal{S}$ be a system of representatives of the isoclasses of
indecomposables of $\T$. We write $\ind\T$ for the full subcategory of
$\T$ whose set of objects is $\mathcal{S}$. For a tree $\Delta$, we write
$\kd$ for the mesh category (see \cite{Rie}). Using the same proof as Riedtmann
\cite{Rie}, one shows the following theorem:

\begin{thm}\label{covering}
There exists a $k$-linear functor $F:\kd\rightarrow \ind\T$
which is surjective and induces  bijections:
$$\bigoplus_{Fz=Fy} \Hom_{\kd}(x,z)\rightarrow \Hom_\T(Fx,Fy),$$
for all vertices $x$ and $y$ of $\mathbb{Z}\Delta$.
\end{thm}  

\subsection{Infinite case}
If the category $\T$ is locally finite not finite, the constructed
functor $F$ is immediately fully faithful. Thus we get the corollary.
\begin{cor}
If $\ind \T$ is not finite, then we have a $k$-linear equivalence
between $\T$ and the mesh category $k(\mathbb{Z}\Delta)$. 
\end{cor}

\subsection{Uniqueness criterion}
The covering functor $F$ can be see as a $k$-linear functor from the derived category $\mathcal{D}^b(\modd k\Delta)$ to the category $\T$. By construction, it satisfies the following property called the \emph{AR-property}:

For each AR-triangle $\xymatrix{X\ar[r]^f & Y\ar[r]^g & Z\ar[r]^h & SX}$ of $\mathcal{D}^b(\modd k\Delta)$, there exists a triangle of $\T$ of the form $\xymatrix{FX\ar[r]^{Ff} & FY\ar[r]^{Fg} & FZ\ar[r]^\epsilon & SFX}$.

In fact, thanks to this property, $F$ is determined by its restriction to the subcategory $\proj k\Delta=k(\Delta)$, \emph{i.e.} we have the following lemma:

\begin{lem}\label{unicity}
Let $F$ and $G$ be $k$-linear functors from $\mathcal{D}^b(\modd k\Delta)$ to $\T$. Suppose that $F$ and $G$ satisfy the AR-property and that the restrictions $F_{|_{k(\Delta)}}$ and $G_{|_{k(\Delta)}}$ are isomorphic. Then the functors $F$ and $G$ are isomorphic as $k$-linear functors. 
\end{lem}
\begin{proof}
It is easy to construct this isomorphism by induction using the (TR3) axiom of the triangulated categories (see \cite{Nee}). 
\end{proof}

\section{Particular cases of $k$-linear equivalence}\label{partcase}

From now we suppose that the category $\T$ is finite, \emph{i.e.} $\T$
has finitely many isoclasses of indecomposable objects.

\subsection{Equivalence criterion}\label{crit}

Let $\Gamma$ be the AR-quiver of $\T$ and suppose that it is isomorphic to
$\mathbb{Z}\Delta/G$. Let $\varphi$ be a generator of $G$. It induces
an automorphism in the mesh category $\kd$ that we still denote by $\varphi$.
Then we have the following equivalence criterion:
\begin{prop}\label{critprop}
The categories $k(\Gamma)$ and $\ind\T$ are equivalent as
$k$-categories if and only if there exists a covering
functor  $F:\kd\rightarrow \ind\T$ and an isomorphism of functors
$\Phi:F\circ\varphi\rightarrow F$. 

\end{prop}
The proof consists in constructing a $k$-linear equivalence between
$\ind\T$ and the orbit category $\kd/\varphi^\mathbb{Z}$ using the
universal property of the orbit category (see
\cite{Kel}), and then constructing an equivalence between
$\kd/\varphi^\mathbb{Z}$ and $k(\Gamma)$.

\subsection{Cylindric case for $\mathbb{A}_n$}

\begin{thm}\label{An}
If $\Delta=\mathbb{A}_n$ and $\varphi=\tau^r$ for some $r\geq 1$, then there exists a
functor isomorphism $\Phi:F\circ\varphi\rightarrow F$, \textit{i.e.}
for each object $x$ of $k(\mathbb{Z}\Delta)$ there exists an automorphism $\Phi_x$ of $Fx$
such that for each arrow $\alpha:x\rightarrow y$ of
$\mathbb{Z}\Delta$, the following diagram commutes:
$$\xymatrix{Fx\ar[r]^{\Phi_x}\ar[d]_{F\alpha} &
  Fx\ar[d]^{F\varphi\alpha}\\ Fy\ar[r]^{\Phi_y} & Fy.}$$
\end{thm} 
To prove this, we need the following lemma:
\begin{lem}
Let $\alpha:x\rightarrow y$ be an arrow of $\mathbb{Z}\mathbb{A}_n$ and let $c$
be a path from $x$ to $\tau^ry$, $r\in\mathbb{Z}$, which is not zero
in the mesh category $k(\mathbb{Z}\mathbb{A}_n)$. 
Then $c$ can be written
$c'\alpha$ where $c'$ is a path from $y$ to $\tau^ry$ (up to sign).
\end{lem}

\begin{proof}
There is a path from $x$ to $\tau^ry$, thus, we have
$\Hom_{\kd}(x,\tau^ry)\simeq k$, and $x$ and $\tau^ry$ are opposite
vertices of a `rectangle' in $\mathbb{Z}\mathbb{A}_n$. This implies that
there exists a path from $x$ to $\tau^r y$ beginning by $\alpha$.
\end{proof}

\begin{proof}\emph{(of theorem \ref{An})} Combining proposition \ref{critprop} and lemma \ref{unicity}, we have just to construct an isomorphism between the restriction of $F$ and $F\circ\varphi$ to a subquiver $\mathbb{A}_n$.

Let us fix a full subquiver of $\mathbb{ZA}_n$ of the following form:
$$\xymatrix{x_1\ar[r]^{\alpha_1} &x_2\ar[r]^{\alpha_2} & \cdots
  \ar[r]^{\alpha_{n-1}} & x_n}$$
such that $x_1\ldots, x_n$ are representatives of the $\tau$-orbits in $\mathbb{Z}\mathbb{A}_n$.
We define the $(\Phi_{x_i})_{i=1\ldots n}$ by induction. We fix
  $\Phi_{x_1}=Id_{Fx_1}$. Now suppose we have constructed some
  automorphisms $\Phi_{x_1},\ldots ,\Phi_{x_i}$ such that for each
  $j\leq i$ the following diagram is commutative:
$$\xymatrix{Fx_{j-1}\ar[r]^{\Phi_{x_{j-1}}}\ar[d]_{F\alpha_{j-1}} &
  Fx_{j-1}\ar[d]^{F\varphi\alpha_{j-1}}\\ Fx_j\ar[r]^{\Phi_{x_j}} &
  Fx_j.}$$
The composition $(F\varphi\alpha_i)\circ\Phi_{x_i}$ is in the morphism
  space $\Hom_\T(Fx_i,Fx_{i+1})$, which is isomorphic, by theorem \ref{covering}, to the space
  $$\bigoplus_{Fz=Fx_{i+1}}\Hom_{\kd}(x_i,z).$$
Thus we can write
  
$$(F\varphi\alpha_i)\Phi_{x_i}=\lambda F\alpha_i + \sum_{z\neq
  x_{i+1}} F\beta_z$$
where $\beta_z$ belongs to $\Hom_{\kd}(x_i,z)$ and $Fz=Fx_{i+1}$.
But $Fz$ is equal to $Fx_{i+1}$ if and only if $z$ is of the form
  $\tau^{rl}x_{i+1}$ for an $l$ in $\mathbb{Z}$. By the lemma, we can
  write $\beta_z=\beta'_z\alpha_i$. Thus we have the equality: 

$$(F\varphi\alpha_i)\Phi_{x_i} = F(\lambda Id_{x_{i+1}} +
\sum_z\beta'_z)F\alpha_i.$$
The scalar $\lambda$ is not zero. Indeed, $\Phi_{x_i}$ is an
automorphism, thus the image of $(F\varphi\alpha_i)\Phi_{x_i}$ is not
zero in the quotient
$$\mathcal{R}(Fx_i,Fx_{i+1})/\mathcal{R}^2(Fx_i,Fx_{i+1}).$$

Thus $\Phi_{x_{i+1}}=F(\lambda Id_{x_{i+1}}+\sum_z\beta'_z)$ is an
automorphism of $Fx_{i+1}$ which verifies the commutation relation
$$(F\varphi\alpha_i)\circ \Phi_{x_i}=\Phi_{x_{i+1}}\circ F\alpha_i.$$

\end{proof}
\subsection{Other standard cases}
In the mesh category $\kd$, where $\Delta$ is a Dynkin tree, the length of
the non zero paths is bounded. Thus there exist automorphisms
$\varphi$ such that, for an arrow $\alpha :x\rightarrow y$ of
$\Delta$, the paths from $x$ to $\varphi^r y$ vanish in the mesh
category for all $r\neq 0$. In other words, for each arrow $\alpha
:x\rightarrow y$ of $\mathbb{Z}\Delta$, we have:
$$\Hom_{\kd/\varphi^{\mathbb{Z}}}(x,y)=\bigoplus_{r\in
  \mathbb{Z}}\Hom_{\kd}(x, \varphi^r y)=\Hom_{\kd}(x,y)\simeq k,$$
where $\kd/\varphi^{\mathbb{Z}}$ is the orbit category (see section
  \ref{crit}). 

\begin{lem}\label{simplelem}
Let $\T$ be a finite triangulated category with AR-quiver
$\Gamma=\mathbb{Z}\Delta/G$. Let $\varphi$ be a generator of $G$ and
suppose that $\varphi$ verifies for each arrow $x\rightarrow y$ of $\mathbb{Z}\Delta$
$$\bigoplus_{r\in
  \mathbb{Z}}\Hom_{\kd}(x, \varphi^r y)=\Hom_{\kd}(x,y)\simeq k.$$
Let $F:\kd\rightarrow \T$ and $G:\kd\rightarrow \T$ be covering functors satisfying the AR-property. Suppose that $F$ and $G$ agree up to isomorphism on the objects of $\kd$. Then $F$ and $G$ are isomorphic as $k$-linear functors.
\end{lem}

\begin{proof}
Using lemma \ref{unicity}, we have just to construct an isomorphism between the functors restricted to $\Delta$.
Let  $\alpha :x\rightarrow y$ be an arrow of $\Delta$. Using
theorem \ref{covering} and the hypothesis, we have the following
isomorphisms:
$$\Hom_\T(Fx,Fy)\simeq \bigoplus_{Fz=Fy}\Hom_{\kd}(x,z) \simeq \bigoplus_{r\in
  \mathbb{Z}}\Hom_{\kd}(x, \varphi^r y)\simeq k$$
and then
$$\Hom_\T(Gx,Gy)\simeq\Hom_\T(Fx,Fy)\simeq k.$$
Thus there exists a scalar $\lambda$ such that
  $G\alpha=\lambda F\alpha$. This scalar does not vanish since $F$ and $G$ are covering functors. As $\Delta$ is a tree, we can find
  some $\lambda_x$ for $x\in \Delta$ by induction such that
$$ G\alpha=\lambda_x\lambda^{-1}_yF\alpha.$$
Now it is easy to check that $\Phi_x=\lambda_xId_{Fx}$ is the functor
isomorphism. 
\end{proof}

This lemma gives us an isomorphism between the functors $F$ and $F\circ \varphi$,  
Moreover, using the same argument, one can show that the covering
functor $F$ is an $S$-functor and a $\tau$-functor.

For each Dynkin tree $\Delta$ we can determine the automorphisms $\varphi$
which satisfy this combinatorial property. Using the preceding lemma
and the equivalence criterion we deduce the following theorem:

\begin{thm}\label{simple}
Let $\T$ be a finite triangulated category with AR-quiver
$\Gamma=\mathbb{Z}\Delta/G$. Let $\varphi$ be a generator of $G$. If one of
these cases holds,
\begin{itemize}
\item
$\Delta=\mathbb{A}_n$ with $n$ odd  and $G$ is generated by $\tau^r$ or
$\varphi=\tau^r\phi$ with $r\geq \frac{n-1}{2}$ and
$\phi=\tau^\frac{n+1}{2}S$;
\item
$\Delta=\mathbb{A}_n$ with $n$ even and $G$ is generated by $\rho^r$
with $r\geq n-1$ and $\rho=\tau^\frac{n}{2}S$;
\item
$\Delta=\mathbb{D}_n$ with $n\geq 5$ and $G$ is generated by $\tau^r$
or $\tau^r\phi$ with $r\geq n-2$ and $\phi$ as in theorem \ref{groupesadmissibles};
\item
$\Delta=\mathbb{D}_4$ and $G$ is generated by $\phi\tau^r$, where $r\geq
2$ and $\phi$ runs over $\sigma_3$;
\item
$\Delta=\mathbb{E}_6$ and $G$ is generated by $\tau^r$ or $\tau^r\phi$
where $r\geq 5$ and $\phi$ is as in theorem \ref{groupesadmissibles};
\item
$\Delta=\mathbb{E}_7$ and $G$ is generated by $\tau^r$, $r\geq 8$;
\item
$\Delta=\mathbb{E}_8$ and $G$ is generated by $\tau^r$, $r\geq 14$.
\end{itemize}
then $T$ is standard, i.e.
the categories $\T$ and $k(\Gamma)$ are equivalent as $k$-linear categories.

\end{thm}

\begin{cor}\label{cycor}
A finite maximal $d$-Calabi-Yau (see \cite[8]{Kel}) triangulated category
$\T$, with $d\geq 2 $, is standard, \emph{i.e.} there exists a $k$-linear
equivalence between $\T$ the orbit category $\mathcal{D}^b(\modd k\Delta)/\tau^{-1}S^{d-1}$ where $\Delta$ is 
Dynkin of type $\mathbb{A}$, $\mathbb{D}$ or $\mathbb{E}$
\end{cor}

\section{Algebraic case}\label{algsection}

For some automorphism groups $G$, we know the $k$-linear structure of
$\T$. But what about the triangulated structure?
We can only give an answer adding hypothesis on the triangulated
structure. In this section, we
distinguish two cases: 

If $\T$ is locally finite, not finite, we have the following theorem
which is proved in section \ref{proof1}:

\begin{thm}\label{algcase}
Let $\T$ be a connected locally finite triangulated category with
infinitely many indecomposables. If $\T$ is the base of a tower of
triangulated categories \cite{Kel3}, then $\T$ is
triangle equivalent to $\mathcal{D}^b(\modd k\Delta)$ for some Dynkin diagram
$\Delta$.
\end{thm}

Now if $\T$ is a finite standard category which is algebraic, \emph{i.e.}  $\T$
is triangle equivalent to $\underline{\mathcal{E}}$ for some
$k$-linear Frobenius category $\mathcal{E}$ (\cite[3.6]{Kel2}), then we have the following result
 which is proved in section \ref{proof2}:

\begin{thm}\label{algcy}
Let $\T$ be a finite triangulated category, which is connected, algebraic and
standard. Then, there exists a Dynkin diagram $\Delta$ of type
$\mathbb{A}$, $\mathbb{D}$ or $\mathbb{E}$ and an auto-equivalence
  $\Phi$ of $\mathcal{D}^b(\modd k\Delta)$ such that $\T$ is triangle equivalent
  to the orbit category $\mathcal{D}^b(\modd k\Delta)/\Phi.$
 
\end{thm}

This theorem combined with corollary \ref{cycor} yields the following
result (compare to \cite[Cor 8.4]{Kel}):

\begin{cor}\label{corcy}
If $\T$ is a finite algebraic maximal $d$-Calabi-Yau category with
$d\geq 2$, then $\T$ is
triangle equivalent to the orbit category
$\mathcal{D}^b(\modd k\Delta)/S^d\nu^{-1}$ for some Dynkin diagram $\Delta$.
\end{cor}

\subsection{$\partial$-functor}
We recall the following definition from \cite{Kel3} and \cite{Ver}.

\begin{df}
Let $\mathcal{H}$ be an exact category and $\T$ a triangulated
category. A \emph{$\partial$-functor} $(I,\partial):
\mathcal{H}\rightarrow \T$ is given by:
\begin{itemize}
\item an additive $k$-linear functor $I:\mathcal{H}\rightarrow \T$;
\item
for each conflation $\xymatrix{\epsilon : X\ \ar@{>->}[r]^i
  & Y\ar@{>>}[r]^p & Z}$ of $\mathcal{H}$, a morphism $\partial \epsilon :IZ\rightarrow SIX$ functorial in $\epsilon$ such that
$\xymatrix{IX\ar[r]^{Ii} & IY\ar[r]^{Ip} &
  IZ\ar[r]^{\partial\epsilon} & SIX}$ is a triangle of $\T$.
\end{itemize} 
\end{df}
For each exact category $\mathcal{H}$, the inclusion
$I:\mathcal{H}\rightarrow \mathcal{D}^b(\mathcal{H})$ can be completed
to a $\partial$-functor $(I,\partial)$ in a unique way.
Let $\T$ and $\T'$ be triangulated categories. If $(F,\varphi):\T\rightarrow \T'$ is an $S$-functor and  $(I,\partial):
\mathcal{H}\rightarrow \T$ is a $\partial$-functor, we say that \emph{$F$ respects
$\partial$} if $(F\circ I,
\varphi(F\partial)):\mathcal{H}\rightarrow \T'$ is a
$\partial$-functor. Obviously each triangle functor respects $\partial$.

\begin{prop}
Let $\mathcal{H}$ be a $k$-linear hereditary abelian category and let $(I,\partial):
\mathcal{H}\rightarrow \T$ be a $\partial$-functor. Then there exists a
unique (up to isomorphism) $k$-linear $S$-functor $F:\mathcal{D}^b(\mathcal{H})\rightarrow
\T$ which respects $\partial$.
\end{prop}    

\begin{proof}
On $\mathcal{H}$ (which can be seen as a full subcategory of
$\mathcal{D}^b(\mathcal{H})$), the functor $F$ is uniquely
determined. We want $F$ to be an $S$-functor, so $F$ is uniquely
determined on $S^n\mathcal{H}$ for $n\in\mathbb{Z}$ too. Since
$\mathcal{H}$ is hereditary, each object of
$\mathcal{D}^b(\mathcal{H})$ is isomorphic to a direct sum of stalk
complexes, \emph{i.e.} complexes concentrated in a single
degree. Thus, the functor $F$ is uniquely determined on the
objects. Now, let $X$ and $Y$ be stalk complexes of
$\mathcal{D}^b(\mathcal{H})$ and $f:X\rightarrow Y$ a non-zero
morphism.
 We can
suppose that $X$ is in $\mathcal{H}$ and $Y$ is in
$S^n\mathcal{H}$. If $n\neq 0,1$, $f$ is necessarily zero. If $n=0$,
then $f$ is a morphism in $\mathcal{H}$ and $Ff$ is uniquely determined. If $n=1$, $f$ is an element of
$\Ext^1_\mathcal{H}(X,S^{-1}Y)$, so gives us a conflation
 $\xymatrix{\epsilon : S^{-1}Y\ \ar@{>->}[r]^i
  & E\ar@{>>}[r]^p & X}$ in $\mathcal{H}$. The functor $F$ respects
 $\partial$, thus $Ff$ has to be equal to $\varphi\circ\partial
\epsilon$  where $\varphi$ is the natural isomorphism between
$SFS^{-1}Y$ and $FY$. Since $\partial$ is functorial, $F$ is a functor.
The result follows.  
\end{proof}

A priori this functor is not a triangle functor.
We recall a theorem proved by B. Keller \cite[cor 2.7]{Kel3}.
\begin{thm}\label{kel} 
Let $\mathcal{H}$ be a $k$-linear exact category, and $\T$ be the base
of a tower of triangulated categories \cite{Kel3}. Let  $(I,\partial):
\mathcal{H}\rightarrow \T$ be a $\partial$-functor such that for each
$n<0$, and all objects
$X$ and $Y$ of $\mathcal{H}$, the space $\Hom_\T(IX,S^nIY)$
vanishes. Then there exists a triangle functor
$F:\mathcal{D}^b(\mathcal{H})\rightarrow \T$ such that the following
diagram commutes up to isomorphism of $\partial$-functors:
$$\xymatrix{\mathcal{H}\ \ar@{^{(}->}[rr] \ar[dr]_{(I,\partial)} & &
  \mathcal{D}^b(\mathcal{H})\ar[dl]^F \\ & \T & }$$ 
\end{thm}

From theorem \ref{kel}, and the proposition above we deduce the
following corollary:
\begin{cor} (compare to \cite{Rin2})
Let $\T$, $\mathcal{H}$ and $(I,\partial):\mathcal{H}\rightarrow \T$
be as in theorem \ref{kel}. If $\mathcal{H}$ is hereditary, then the unique functor $F:\mathcal{D}^b(\mathcal{H})\rightarrow \T$
which respects $\partial$ is a triangle functor.
\end{cor}

\subsection{Proof of theorem \ref{algcase}}\label{proof1}
Let $F$ be the $k$-linear equivalence constructed in theorem \ref{covering} between an algebraic triangulated
category $\T$ and $\mathcal{D}^b(\mathcal{H})$ where
$\mathcal{H}=\modd k\Delta$ and $\Delta$ is a simply laced Dynkin
graph. As we saw in section \ref{partcase}, the covering functor is an $S$-functor.

The category $\mathcal{H}$ is the heart of the standard $t$-structure
on $\mathcal{D}^b(\mathcal{H})$. The image of this $t$-structure
through $F$ is a $t$-structure on $\T$. Indeed, $F$ is an $S$-equivalence, so
the conditions $(i)$ and $(ii)$ from \cite[Def 1.3.1]{Bei} hold obviously. And since
$\mathcal{H}$ is hereditary, for an object $X$ of
$\mathcal{D}^b(\mathcal{H})$, the morphism $X\rightarrow S\tau_{\leq
  0}X$ of the triangle
$$\xymatrix{\tau_{\leq0}X\ar[r]  & X\ar[r] & \tau_{\geq 0} X \ar[r]&
S\tau_{\leq 0} X}$$ vanishes. Thus the image of this triangle  through
$F$ is a triangle of $\T$ and condition $(iii)$ of \cite[Def 1.3.1]{Bei} holds. Then we get a $t$-structure on $\T$ whose
heart is $\mathcal{H}$.

It results from \cite[Prop 1.2.4]{Bei} that the inclusion of the
heart of a $t$-structure can be uniquely completed to a
$\partial$-functor. Thus we obtain a $\partial$-functor
$(F_0,\partial):\mathcal{H}\rightarrow \T$ with $F_0= F_{|_{\mathcal{H}}}$.

 The functor $F$ is an $S$-equivalence. Thus for each $n<0$, and all
  objects $X$
  and $Y$ of $\mathcal{H}$, the space $\Hom_\T(FX,S^nFY)$ vanishes. Now we can apply theorem \ref{kel} and we get the following
  commutative diagram:
$$\xymatrix{\mathcal{H}\ \ar@{^{(}->}[rr] \ar[dr]_{(F_0,\partial)} & &
  \mathcal{D}^b(\mathcal{H}),\ar@<.5ex>[dl]^F \ar@<-.5ex>[dl]_{G}\\ &
  \T & }$$
where $F$ is the $S$-equivalence and $G$ is a triangle functor. Note
  that a priori $F$ is an $S$-functor which does not respect $\partial$.
The functors $F_{|_\mathcal{H}}$ and $G_{|_\mathcal{H}}$ are
  isomorphic. The functor $F$ is an $S$-functor thus we have an
  isomorphism $F_{|_{S^n\mathcal{H}}}\simeq G_{|_{S^n\mathcal{H}}}$
  for each $n\in \mathbb{Z}$. Thus the functor $G$ is essentially
  surjective. Since $\mathcal{H}$ is the category $\modd k\Delta$, to
  show that $G$ is fully faithful, we have
  just to show that for each $p\in\mathbb{Z}$, there is an isomorphism
  induced by $G$
  $$\xymatrix{\Hom_{\mathcal{D}^b(\mathcal{H})}(A,S^pA)\ar[r]^\sim &
  \Hom_{\T}(GA,S^pGA)}$$ where $A$ is the free module $k\Delta$. For
  $p=0$, this is clear because $A$ is in $\mathcal{H}$. And for $p\neq
  0$ both sides vanish. 

Thus $G$ is a triangle equivalence between
$\mathcal{D}^b(\mathcal{H})$ and $\T$.

\subsection{Finite algebraic standard case}\label{proof2}
For a small dg category $\ca$, we denote by $\cc\ca$ the category of dg $\ca$-modules,
by $\cd\ca$ the derived category of $\ca$ and by $\per\ca$ the
{\em perfect derived category} of $\ca$, \emph{i.e.} the smallest triangulated subcategory
of $\cd\ca$ which is stable under passage to direct factors and contains
the free $\ca$-modules $\ca(?,A)$, where $A$ runs through the objects of $\ca$. Recall that a
small triangulated category is {\em algebraic} if it is triangle
equivalent to $\per\ca $ for a dg category $\ca$. For two small dg categories
$\ca$ and $\cb$, a triangle functor $\per\ca\to\per\cb$ is {\em algebraic} if
it is isomorphic to the functor
\[
F_X=?\lten_\ca X
\]
associated with a dg bimodule $X$, \emph{i.e.} an object of the derived category $\cd(\ca^{op}\ten\cb)$.

Let $\Phi$ be an algebraic autoequivalence of $\cd^b(\modd k\Delta)$ such
that the orbit category $\cd^b(\modd k\Delta)/\Phi$ is triangulated.
Let $Y$ be a dg $k\Delta$-$k\Delta$-bimodule such that $\Phi=F_Y$. In
section~9.3 of \cite{Kel}, it was shown that there is a canonical triangle
equivalence between this orbit category and the perfect derived category
of a certain small dg category.
Thus, the orbit category is algebraic, and endowed with a canonical triangle
equivalence to the perfect derived category of a small dg category. Moreover,
by the construction in [loc. cit.], the projection functor
\[
\pi : \cd^b(\modd k\Delta) \to \cd^b(\modd k\Delta)/\Phi
\]
is algebraic.

The proof of theorem 7.0.5 is based on the following universal property of
the triangulated orbit category $\cd^b(\modd k\Delta)/\Phi$. For the
proof, we refer to section~9.3 of \cite{Kel}.

\begin{prop} Let $\cb$ be a small dg category and
\[
F_X=?\lten_{k\Delta} X: \cd^b(\modd k\Delta) \to \per\cb
\]
an algebraic triangle functor given by a dg $k\Delta$-$\ca$-bimodule $X$.
Suppose that there is an isomorphism between $Y\lten_{k\Delta} X$ and
$X$ in the derived bimodule category $\cd(k\Delta^{op}\ten\cb)$.
Then the functor $F_X$ factors, up to isomorphism of triangle functors,
through the projection
\[
\pi : \cd^b(\modd k\Delta) \to \cd^b(\modd k\Delta)/\Phi .
\]
Moreover, the induced triangle functor is algebraic.
\end{prop}

Let us recall a lemma of Van den Bergh \cite{Kel4}:
\begin{lem}\label{vdb}
Let $Q$ be a quiver without oriented cycles and $\mathcal{A}$ be a dg category. We denote by $k(Q)$ the category of paths of $Q$ and by $Can:\mathcal{CA}\rightarrow \mathcal{DA}$ the canonical functor. Then we have the following properties:

\hspace{.2in} $a)$  Each functor $F:k(Q)\rightarrow \mathcal{DA}$ lifts, up to isomorphism, to a functor $\tilde{F}:k(Q)\rightarrow \mathcal{CA}$ which verifies the following property: For each vertex $j$ of $Q$, the induced morphism
$$\bigoplus_{i}\tilde{F}i\rightarrow \tilde{F}j,$$
where $i$ runs through the immediate predecessors of $j$, is a monomorphism which splits as a morphism of graded $\mathcal{A}$-modules.

\hspace{.2in} $b)$
Let $F$ and $G$ be functors from $k(Q)$ to $\mathcal{CA}$, and suppose that $F$ satisfies the property of $a)$. Then any morphism of functors $\varphi:Can\circ F\rightarrow Can \circ G$ lifts to a morphism $\tilde{\varphi}:F\rightarrow G$.

\end{lem}

\begin{proof}
 $a)$
 For each vertex $i$ of $Q$, the object $Fi$  is isomorphic in $\mathcal{DA}$ to its cofibrant resolution $X_i$. Thus for each arrow $\alpha: i\rightarrow j$, $F$ induces a morphism $f_\alpha:X_i\rightarrow X_j$ which can be lifted to $\mathcal{CA}$ since the $X_i$ are cofibrant. Since $Q$ has no oriented cycle, it is easy to choose the $f_\alpha$ such that the property is satisfied.

\hspace{.2in} $b)$ For each vertex $i$ of $Q$, we may assume that $Fi$ is cofibrant. Then we can lift $\varphi_i:Can \circ Fi\rightarrow Can \circ Gi$ to $\psi_i:Fi\rightarrow Gi$. For each arrow $\alpha$ of $Q$, the square
$$\xymatrix{Fi\ar[r]^{F_\alpha}\ar[d]^{\psi_i} & Fj\ar[d]^{\psi_j} \\ Gi\ar[r]^{G_\alpha} & Gj}$$
is commutative in $\mathcal{DA}$.
Thus the square $$\xymatrix{\bigoplus_i Fi\ar[r]^{F_\alpha}\ar[d]^{(\psi_i)} & Fj\ar[d]^{\psi_j} \\ \bigoplus_i Gi\ar[r]^{G_\alpha} & Gj}$$
is commutative up to nullhomotopic morphism $h:\bigoplus_i Fi \rightarrow Gj$. Since the morphism $f:\bigoplus_i Fi \rightarrow Fj$ is split mono in the category of graded $\mathcal{A}$-modules, $h$ extends along $f$ and we can modify $\Psi_j$ so that the square becomes commutative in $\mathcal{CA}$. The quiver $Q$ does not have oriented cycles, so we can construct $\tilde{\varphi}$ by induction.
\end{proof}

\begin{proof}\emph{(of theorem \ref{algcy})}
The category $\T$ is small and algebraic, thus we may assume that
$\T=\per \mathcal{A}$ for some small dg category $\mathcal{A}$. Let
$F:\mathcal{D}^b(\modd k\Delta)\rightarrow \T$ be the
covering functor of theorem \ref{covering}. Let  $\Phi$ be an auto-equivalence of $\mathcal{D}^b(\modd k\Delta)$ such
that the AR-quiver of the orbit category $\mathcal{D}^b(\modd k\Delta)/\Phi$ is isomorphic (as translation quiver) to the AR-quiver of $\T$. We may
assume that $\Phi=-\lten_{k\Delta}Y$ for an object $Y$ of
$\mathcal{D}(k\Delta^{op}\ten k\Delta)$. The orbit category $\mathcal{D}^b(\modd k\Delta)/\Phi$ is algebraic, thus it is $\per\mathcal{B}$ for some dg category
$\mathcal{B}$.

The functor $F_{|_{k(\Delta)}}$ lifts by lemma \ref{vdb} to a functor $\tilde{F}$ from $k(\Delta )$ to $\mathcal{CA}$. This means that the object $X=\tilde{F}(k\Delta)$ has a structure of dg $k\Delta^{op}\ten \mathcal{A}$-module. We denote by $X$ the image of this object in $ \mathcal{D}(k\Delta^{op}\ten \mathcal{A})$. 

The functors $F$ and $-\lten_{k\Delta}X$ become isomorphic when restricted to $k(\Delta)$. Moreover  $-\lten_{k\Delta}X$ satisfies the AR-property since it is a triangulated functor. Thus by lemma \ref{unicity}, they are isomorphic as $k$-linear functors. So we have the following diagram:
$$ \xymatrix{ & \mathcal{D}^b(\modd k\Delta) 
  \ar[rr]^{-\lten_{k\Delta}X} \ar@(dl,dr)_{-\lten_{k\Delta}Y} & & \per\mathcal{A}=\T \\ &
  & & } $$

The category $\T$ is standard, thus there exists an isomorphism of $k$-linear functors:
$$\xymatrix{c:-\lten_{k\Delta}X\ar[rr] & &
-\lten_{k\Delta}Y \lten_{k\Delta}X.}$$
The functor $-\lten_{k\Delta}X$ restricted to the category $k(\Delta)$ satisfies the property of $a)$ of lemma \ref{vdb}. Thus we can apply $b)$ and lift $c_{|_{k(\Delta)}}$ to an isomorphism $\tilde{c}$ between $X$ and $Y\lten_{k\Delta}X$ as dg-$k\Delta^{op}\ten \mathcal{A}$-modules.

By the universal property of the orbit category, the bimodule $X$ endowed
with the isomorphism $\tilde{c}$ yields a triangle functor from
$\mathcal{D}^b(\modd k\Delta)/\Phi$ to $\T$ which comes from a
bimodule $Z$ in $\mathcal{D}(\mathcal{B}^{op}\ten
\mathcal A)$.
 
$$ \xymatrix{  \mathcal{D}^b(\modd k\Delta) 
  \ar[rr]^{-\lten_{k\Delta}X}\ar[d]_\pi
  \ar@(ul,ur)^{-\lten_{k\Delta}Y} & & \per\mathcal{A}=\T \\
  \mathcal{D}^b(\modd k\Delta)/\Phi=\per\mathcal{B}\ar@{-->}[urr]_{-\lten_{k\Delta}Z} &
   & } $$
The functor $-\lten_{k\Delta}Z$ is essentially surjective. Let us show that it is fully faithful. For $M$ and $N$ objects of $\mathcal{D}^b(\modd k\Delta)$ we have the following commutative diagram:
$$\xymatrix{& \bigoplus_{n\in \mathbb{Z}} \Hom_\mathcal{D} (M,\Phi^nN)\ar[dr]^{-\lten_{k\Delta}X=F} \ar[dl]_\pi& \\ \Hom_{\mathcal{D}/\Phi}(\pi M,\pi N)\ar[rr]^{-\lten_{k\Delta}Z} & & \Hom_\T(FM,FN),}$$
where $\mathcal{D}$ means $\mathcal{D}^b(\modd k\Delta)$. The two diagonal morphisms are isomorphisms, thus so is the horizontal morphism. This proves that $-\lten_{k\Delta}Z$ is a triangle equivalence between the orbit category $\mathcal{D}^b(\modd k\Delta)/\Phi$ and $\T$.
\end{proof}

\section{Triangulated structure on the category of projectives}\label{proj}

Let $k$ be a algebraically closed field and $\cp$ a $k$-linear category
with split idempotents. The
category $\modd\cp$ of contravariant finitely
presented functors from
$\cp$ to $\modd k$ is exact. As the idempotents split, the projectives
of $\modd\cp$ coincide with the representables. Thus the Yoneda functor
gives a natural equivalence between $\cp$ and $\proj \cp$. Assume
besides that $\modd \cp$ has a structure of Frobenius category. The stable category
$\underline{\modd} \cp$ is a triangulated category, we write $\Sigma$
for the suspension functor.

Let $S$ be an auto-equivalence of $\cp$. It can be extended to an exact
functor from $\modd\cp$ to $\modd\cp$ and thus to a triangle functor
of $\underline{\modd}\cp$. The aim of this part is to find a necessary
condition on the functor $S$ such that the category $(\cp, S)$ has a
triangulated structure. Heller already showed \cite[thm 16.4]{Hel} that
if there exists an isomorphism of triangle functors between $S$ and
$\Sigma^3$, then $\cp$ has a pretriangulated structure. But he did not
succeed in proving the octahedral axiom. We are going to impose a stronger
condition on the functor $S$ and prove the following theorem:

\begin{thm}\label{thmprincipal}
Assume there exists an exact sequence of exact functors from $\modd\cp$
to $\modd\cp$:  
$$\xymatrix{0\ar[r] & Id\ar[r] & X^0\ar[r] & X^1\ar[r] &
  X^2\ar[r] & S\ar[r] &0},$$
where the $X^i$, $i=0,1,2$,  take values in $\sf{proj}\cp$. Then the category $\cp$
  has a structure of triangulated category with suspension functor
  $S$.
\end{thm}
For an $M$ in $\modd\cp$, denote  $\xymatrix{T_M: X^0M\ar[r]& X^1M\ar[r] & X^2M\ar[r] &
  SX^0M}$ a \emph{standard triangle}. A triangle of $\cp$
  will be a sequence $\xymatrix{X:P\ar[r]^u & Q\ar[r]^v & R\ar[r]^w &
  SP}$ which is isomorphic to a standard triangle $T_M$ for an $M$ in
  $\modd\cp$.

\subsection{S-complexes, $\Phi$-S-complexes and standard triangles}

Let $\mathcal{A}cp(\modd\cp)$ be the category of acyclic complexes with
projective components. It is a Frobenius category whose 
projective-injectives are the contractible complexes, \emph{i.e.} the
complexes homotopic to zero. The functor
$Z^0:\mathcal{A}cp(\modd\cp)\rightarrow \modd\cp$ which sends a complex
$$\xymatrix{ \cdots \ar[r] & X^{-1}\ar[r]^{x^{-1}} & X^0\ar[r]^{x^0} &
  X^1\ar[r]^{x^1} &\cdots}$$ to the kernel of $x^0$ is an exact
  functor. It sends the projective-injectives to
  projective-injectives and induces a triangle equivalence between
  $\underline{\mathcal{A}cp(\modd\cp)}$ and $\underline{\modd}\cp$.  
\begin{df}
An object of $\mathcal{A}cp(\modd\cp)$ is called an \emph{$S$-complex} if
it is $S$-periodic, \emph{i.e.} if
it has the following form:
$$\xymatrix{\cdots \ar[r] &P\ar[r]^u & Q\ar[r]^v & R\ar[r]^w &
  SP\ar[r]^{Su} & SQ\ar[r] &\cdots}.$$
\end{df}
The category $S$-\comp of $S$-complexes with $S$-periodic morphisms is
a non full subcategory of $\mathcal{A}cp(\modd\cp)$. It is a Frobenius
category. The projective-injectives are the $S$-contractibles,
\emph{i.e.} the complexes homotopic to zero with an $S$-periodic
homotopy. Using the functor $Z^0$, we get an exact functor from $S$-\comp to
$\modd\cp$ which induces a triangle functor:  
$$\underline{Z^0}: \underline{S\textrm{-\comp}}\longrightarrow \underline{\modd}\cp.$$

Fix a sequence  as in theorem \ref{thmprincipal}. Clearly, it induces for each object $M$ of $\underline{\modd}\cp$, a functorial
  isomorphism in $\underline{\modd}\cp$, 
$\xymatrix{\Phi_M:\Sigma^3 M\ar[r] & SM}.$

Let $Y$ be an $S$-complex, $$\xymatrix{Y: \cdots \ar[r] &P\ar[r]^u & Q\ar[r]^v & R\ar[r]^w &
  SP\ar[r]^{Su} & SQ\ar[r] &\cdots}.$$ Let $M$ be the kernel of $u$.
  Then $Y$ induces an isomorphism $\theta$ (in $\underline{\modd}\cp$) between $\Sigma^3 M$ and
  $SM$. If $\theta$ is equal to $\Phi_M$, we will say that $X$ is a \emph{$\Phi$-$S$-complex}. 

Let $M$ be an object of $\modd\cp$. The standard triangle $T_M$ can be
  see as a $\Phi$-$S$-complex: $$\xymatrix{\cdots\ar[r] & X^0M\ar[r] &
  X^1M\ar[r] & X^2M\ar[r] & SX^0M\ar[r] & SX^1M\ar[r] & \cdots.}$$

The functor $T$ which sends an object $M$ of $\modd\cp$ to the
$S$-complex $T_M$ is exact since the $X^i$ are exact. It satisfies the
relation $Z^0\circ T\simeq Id_{\modd\cp}.$ Moreover, as it preserves
the projective-injectives, it induces a triangle functor:
$$T:\underline{\modd}\cp\rightarrow 
  \underline{S\textrm{-\comp}}.$$

\subsection{Properties of the functors $Z^0$ and $T$}
\begin{lem}\label{lem0}
An $S$-complex which is homotopy-equivalent to a $\Phi$-$S$-complex is
a $\Phi$-$S$-complex.
\end{lem}
\begin{proof}
Let $\xymatrix@1{X:P\ar[r]^(.6)u & Q\ar[r]^v & R\ar[r]^w & SP}$ be an $S$-complex
homotopy-equivalent to the $\Phi$-$S$-complex
$\xymatrix@1{X':P'\ar[r]^(.6){u'} & Q'\ar[r]^{v'} & R'\ar[r]^{w'} & SP'}$. Let
$M$ be the kernel of $u$ and $M'$ the kernel of $u'$. By assumption, there
exists a $S$-periodic homotopy equivalence $f$ from $X$ to $X'$, which induces a
morphism $g=Z^0f:M\rightarrow M'$. Thus, we get the following
commutative diagram:
$$\xymatrix{ & P\ar[r]\ar[dd]^{f^0} & Q\ar[r]\ar[dd]^{f^1} &
  R\ar[rrr]\ar[dd]^{f^2}\ar@{>>}[dr] & & & SP\ar[dd]^{Sf^0}\\
M\quad \ar@<-1ex>[dd]^g\ar@{>->}[ur] & & & & \Sigma^3M\ar[r]^{\theta}\ar[dd]^(.3){\Sigma^3g} &
  SM\quad \ar@{>->}[ur]\ar@<-1ex>[dd]^(.3){Sg} & \\
& P'\ar[r] & Q'\ar[r] & R'\ar@{>>}[dr]\ar@{-}[r] & \ar@{-}[r] & \ar[r] & SP'.\\
M'\quad \ar@{>->}[ur] & & & & \Sigma^3M'\ar[r]_{\Phi_{M'}} & SM'\quad\ar@{>->}[ur] &}$$ 
The morphism $g$ is an isomorphism of $\underline{\modd}\cp$ since
  $f$ is an isomorphism of $\underline{S\textrm{-\comp}}$. Thus the
  morphisms $\Sigma^3g$ and $Sg$ are isomorphisms of
  $\underline{\modd}\cp$. The following equality in $\underline{\modd}\cp$
$$\theta=(Sg)^{-1}\Phi_{M'}\Sigma^3g=\Phi_M$$ shows that the complex
$X$ is a $\Phi$-$S$-complex.
\end{proof}

\begin{lem}\label{lem1}
Let $$\xymatrix{X:P\ar[r]^u & Q\ar[r]^v & R\ar[r]^w & SP} \textrm{
  and }
\xymatrix{X':P'\ar[r]^{u'} & Q'\ar[r]^{v'} & R'\ar[r]^{w'} & SP'}$$
be two $\Phi$-$S$-complexes. Suppose that we have a commutative square:
$$\xymatrix{P\ar[r]^u\ar[d]_{f^0} & Q\ar[d]^{f^1} \\ P'\ar[r]^{u'} &
  Q'.}$$
Then, there exists a morphism $f^2:R\rightarrow R'$ such that
  $(f^0,f^1,f^2)$ extends to an $S$-periodic morphism from $X$ to $X'$.
\end{lem}

\begin{proof}
Let $M$ be the kernel of $u$, $M'$ be the kernel of $u'$ and
$f:M\rightarrow M'$ be the morphism induced by the commutative
square. As $R$ and $R'$ are projective-injective objects, we can find
a morphism $g^2:R\rightarrow R'$ such that the following square commutes: 
$$\xymatrix{Q\ar[r]^v\ar[d]_{f^1} & R\ar[d]^{g^2}\\ Q'\ar[r]^{v'} &
  R'.}$$
The morphism $g^2$ induces a morphism $g:SM\rightarrow SM'$ such that
  the following square is commutative in $\underline{\modd}\cp$:
$$\xymatrix{\Sigma^3M\ar[r]^{\Phi_M}\ar[d]_{\Sigma^3f} & SM\ar[d]^g\\ \Sigma^3M'\ar[r]^{\Phi_{M'}} &
  SM'.}$$
Thus the morphisms $Sf$ and $g$ are equal in $\underline{\modd}\cp$,
  \emph{i.e.} there exists a projective-injective $I$ of $\modd\cp$ and
  morphisms $\alpha:SM\rightarrow I$ and $\beta:I\rightarrow SM'$ such
  that $g-Sf=\beta\alpha$. Let $p$ (resp. $p'$) be the epimorphism
  from $R$ onto $SM$ (resp. from $R'$ onto $SM'$). Then, as $I$ is
  projective, $\beta$ factors through $p'$. 
$$\xymatrix{Q\ar[r]^v\ar[ddd]_{f^1} & R\ar[rr]^w\ar@{>>}[dr]_p
  \ar[ddd]_{g^2} & & SP\ar[ddd]^{Sf^0}\\
& & SM\quad\ar@{>->}[ur]\ar[d]_\alpha \ar@(dr,ur)[ddd]^{Sf} \ar@(dl,ul)[ddd]_(.3){g}& \\
& & I\ar[dd]_(.3)\beta \ar@{.>}[dl]_\gamma & \\
Q'\ar[r]^{v'} & R'\ar@{-}[r]\ar@{>>}[dr]_{p'} & \ar[r] & SP'\\
& & SM'\quad\ar@{>->}[ur] & }$$

We put $f^2=g^2-\gamma\alpha p$. Then obviously, we have the
  equalities $f^2v=v'f^1$
  and $w'f^2=Sf^0w$. Thus the morphism $(f^0,f^1,f^2)$ extends to a
  morphism of $S$-$\comp$.
\end{proof} 

\begin{prop}\label{thm8}
The functor
$\underline{Z^0}:\underline{\Phi\textrm{-S-\emph{\comp}}}\longrightarrow\underline{\modd}\cp$
is full and essentially surjective. Its kernel is an ideal whose
square vanishes.
\end{prop}
\begin{proof}
The functor $\underline{Z^0}$ is essentially surjective since we have the
relation $\underline{Z^0}\circ \underline{T}=Id_{\underline{\modd}\cp}.$

Let us show that $\underline{Z^0}$ is full. Let
 $$\xymatrix{X:P\ar[r]^u & Q\ar[r]^v & R\ar[r]^w & SP} \textrm{
  and }
\xymatrix{X':P'\ar[r]^{u'} & Q'\ar[r]^{v'} & R'\ar[r]^{w'} & SP'}$$
be two $\Phi$-$S$-complexes. Let $M$ (resp. $M'$) be the kernel of $u$
  (resp. $u'$). As $P$, $Q$, $P'$ and $Q'$ are projective-injective,
  there exist morphisms $f^0:P\rightarrow P'$ and $f^1:Q\rightarrow
  Q'$ such that the following diagram commutes:
$$\xymatrix{M\ \ar@{>->}[r]\ar[d]_f & P\ar[r]^u\ar[d]^{f^0} &
  Q\ar[d]^{f^1}\\ M'\ \ar@{>->}[r] & P'\ar[r]^{u'} & Q'.}$$
Now the result follows from lemma \ref{lem1}.

Now let $\underline{f}:X\rightarrow X'$ be a morphism in the kernel of
$\underline{Z^0}$. Up to homotopy, we can suppose that $\underline{f}$
has the following form:
$$\xymatrix{P\ar[r]^u\ar[d]_0 & Q\ar[r]^v\ar[d]_0 & R\ar[r]^w\ar[d]_{f^2}
  & SP\ar[d]^0\\
P'\ar[r]^{u'} & Q'\ar[r]^{v'} & R'\ar[r]^{w'} & SP'.}$$

As the composition $w'f^2$ vanishes and as $Q'$ is
projective-injective, $f^2$ factors through $v'$. For the same
argument, $f^2$ factors through $w$. If $\underline{f}$ and
$\underline{f'}$ are composable morphisms of the kernel of
$\underline{Z^0}$, we get the following diagram: 

$$\xymatrix{P\ar[r]^u\ar[d]_0 & Q\ar[r]^v\ar[d]_0 & R\ar[r]^w\ar[d]_{f^2}\ar@{.>}[dl]_{h^2}
  & SP\ar[d]^0\\
P'\ar[r]^{u'}\ar[d]_0 & Q'\ar[d]_0\ar[r]^{v'} &
  R'\ar[d]_{f'^2}\ar[r]^{w'} & SP'\ar[d]^0\ar@{.>}[dl]_{h'^3}\\
P''\ar[r]^{u''} & Q''\ar[r]^{v''} & R''\ar[r]^{w''} & SP''.}$$
The composition $\underline{f'}\underline{f}$ vanishes obviously.

\end{proof} 

\begin{cor}\label{cor}
A $\Phi$-$S$-complex morphism $f$ which induces an isomorphism
$\underline{Z^0}(f)$ in $\underline{\modd}\cp$ is an homotopy-equivalence.  
\end{cor}

This corollary comes from the previous theorem and from the following lemma.
 
\begin{lem}
Let $F:\mathcal{C}\rightarrow\mathcal{C'}$ be a full functor between two additive categories. If the kernel of
$F$ is an ideal whose square vanishes, then $F$ detects isomorphisms.
\end{lem}

\begin{proof}
Let $u\in\Hom_{\mathcal{C}}(A,B)$ be a morphism in $\mathcal{C}$ such
that $Fu$ is an isomorphism. Since the functor $F$ is full, there exists $v$ in
$\Hom_\mathcal{C}(B,A)$ such that $Fv=(Fu)^{-1}$. The morphism $w=uv-Id_B$
is in the kernel of $F$, thus $w^2$ vanishes. Then the morphism
$v(Id_B-w)$ is a right inverse of $u$. In the same way we show that
$u$ has a left inverse, so $u$ is an isomorphism.   
\end{proof}

\begin{prop}\label{thm11}
The category of $\Phi$-$S$-complexes is equivalent to the category of
$S$-complexes which are homotopy-equivalent to standard triangles.
\end{prop}

\begin{proof}
Since standard triangles are $\phi$-$S$-complexes, each $S$-complex that
is homotopy equivalent to a standard triangle is a $\Phi$-$S$-complex
(lemma \ref{lem0}).

Let $\xymatrix{X:P\ar[r]^u & Q\ar[r]^v & R\ar[r]^w & SP}$ be a 
$\Phi$-$S$-complex. Let $M$ be the kernel of $u$. Then there exist
morphisms $f^1:P\rightarrow X^0M$ and $f^1: Q\rightarrow X^1M$ such
that the following diagram is commutative:
$$\xymatrix{M\ \ar@{>->}[r]\ar@{=}[d] & P\ar[r]^u\ar[d]^{f^0} &
  Q\ar[d]^{f^1}\\ M\ \ar@{>->}[r] & X^0M\ar[r] & X^1M.}$$
We can complete (lemma \ref{lem1}) $f$ into an $S$-periodic morphism from
  $X$ in $T_M$. The morphism $f$ satisfies $Z^0f=Id_M$, so
  $\underline{Z}^0(T_M)$ and $Z^0(X)$ are equal in
  $\underline{\modd}\cp$. By the corollary, $T_M$ and $X$ are
  homotopy-equivalent. Thus the inclusion functor $T$ is essentially surjective.

\end{proof}

These two diagrams summarize the results of this section:

$$\xymatrix{\Phi\textrm{-S-\comp}\ar@{^(->}[rr]^{\textrm{\tiny{full}}} & &
  \textrm{\begin{minipage}{1.4cm}S-\comp
  \tiny{(Frobenius)}\end{minipage}}\ar[rr]\ar@<1ex>[dd]^(.4){
  \textrm{\begin{minipage}{.8cm} \center{$Z^0$}\\\tiny{(exact)}\end{minipage}}}&
  &\textrm{\begin{minipage}{2cm}$\mathcal{A}cp(\modd\cp)$
  \tiny{(Frobenius)}\end{minipage}}\ar[ddll]^{
  \textrm{\begin{minipage}{1.5cm}\center{$Z^0$} \\ \tiny{(exact)}\end{minipage}}}\\
&&&&\\
&&\textrm{ \begin{minipage}{1.3cm}\center{$\modd\cp$}\\ \tiny{(Frobenius)}\end{minipage}}\ar@<1ex>[uu]^(.6){\textrm{\begin{minipage}{.8cm}\center{$T$}\\\tiny{(exact)}\end{minipage}}}\ar[uull]^T
&& }$$

$$\xymatrix{\underline{\Phi\textrm{-S-\comp}} \ar@{^(->}[rr] & &
  \textrm{\begin{minipage}{1.1cm}\underline{S-\comp} \tiny{(triang.)}\end{minipage}}\ar[rr]\ar@<1ex>[dd]^(.3){\textrm{\begin{minipage}{1.7cm}\center{$\underline{Z^0}$}\\
   \tiny{(exact, full, ess. surj., pr.\ref{thm8})}\end{minipage}}}&
  &\textrm{\begin{minipage}{2cm}\underline{$\mathcal{A}cp(\modd\cp)$}\\\tiny{(triang.)}\end{minipage}}\ar[ddll]^(.6){
  \textrm{\begin{minipage}{2cm}\center{$\underline{Z^0}$}\\ \tiny{
  (triangle equivalence)}\end{minipage}}}\\
&&&&\\
&&\textrm{\begin{minipage}{1.1cm}$\underline{\modd}\cp$\\\tiny{(triang.)}\end{minipage}}\ar@<1ex>[uu]^{\textrm{\begin{minipage}{.9cm}\center{$\underline{T}$}\\\tiny{(triang.)}\end{minipage}}}\ar[uull]^{
\textrm{\begin{minipage}{1.2cm}\center{$\underline{T}$}\\\tiny{(ess.surj.,
  pr. \ref{thm11})}\end{minipage}}}&& }$$

\subsection{Proof of theorem \ref{thmprincipal}}

We are going to show that the $\Phi$-$S$-complexes form a system of
triangles of the category $\cp$. We use triangle axioms as in \cite{Nee}.
\\

\textbf{TR0}:
For each object $M$ of $\cp$, the $S$-complex
 $\xymatrix{M\ar@{=}[r] &
  M\ar[r] & 0\ar[r] & SM}$ is homotopy-equivalent to the zero
  complex, so is a $\Phi$-$S$-complex.
\\ 

\textbf{TR1}:
Let $u:P\rightarrow Q$ be a morphism of $\cp$, and let $M$ be its
kernel. We can find morphisms $f^0$ and $f^1$ so as to obtain a
commutative square:
$$
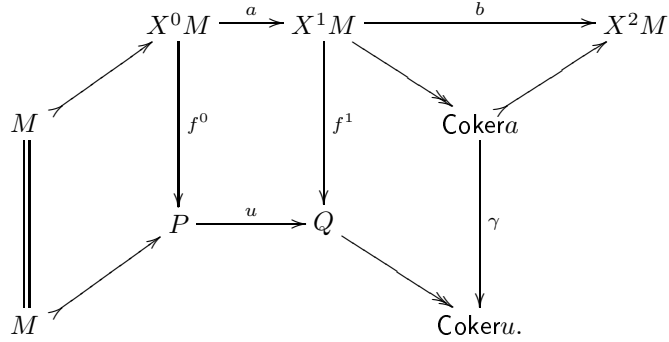
\xymatrix{& X^0M\ar[r]^a \ar[dd]^{f^0} &
  X^1M\ar[rr]^b\ar@{>>}[rd]\ar[dd]^{f^1} & & X^2M\\
M\quad\ar@{>->}[ur]\ar@{=}[dd]<-1ex> & & & \text{\sf{Coker}}a\ar@{>->}[ur]\ar[dd]^{\gamma} & \\
& P\ar[r]^u & Q\ar@{>>}[dr] & & \\ M\quad \ar@{>->}[ur] & & & \text{\sf{Coker}}u.
  & }$$
We form the following push-out:
$$\xymatrix{\ar@{}[drrr]|{\text{PO}} 0\ar[r] & \text{\sf{Coker}}a\ar[r]\ar[d]^{\gamma} & X^2M\ar[d]\ar[r] &
  SM\ar@{=}[d]\ar[r] & 0\\
0\ar[r] & \text{\sf{Coker}}u\ar[r] & R\ar[r] & SM\ar[r] & 0.}$$
It induces a triangle morphism of the triangular category $\underline{\modd}\cp$:
$$\xymatrix{\text{{\sf Coker }}a\ar[r]\ar[d]^{\gamma} & X^2M\ar[r]\ar[d] &
  SM\ar[r]\ar@{=}[d] &\Sigma\text{{\sf Coker }}a\ar[d]^{\Sigma\gamma}\\ \text{{\sf Coker }}u\ar[r] & R\ar[r] &
  SM\ar[r] &\Sigma\text{{\sf Coker }}u.}$$
The morphism $\gamma$ is an isomorphism in $\underline{\modd}\cp$ since
  {\sf Coker }$a$ and {\sf Coker }$u$ are canonically isomorphic to
  $\Sigma^2M$ in $\underline{\modd}\cp$. By the five lemma,
  $X^2M\rightarrow R$ is an isomorphism in
  $\underline{\modd}\cp$. Since $X^2M$ is projective-injective, so is
  $R$. Thus the complex $\xymatrix{P\ar[r]^u&Q\ar[r]&R\ar[r]&SP}$
is an $S$-complex. Then we have to see that it is a $\Phi$-$S$-complex.
Let $\theta$ be the isomorphism between $SM$ and $\Sigma^3M$ induced
  by this complex. We write $\alpha$ (resp. $\beta$) for the canonical
  isomorphism in $\underline{\modd}\cp$ between $\Sigma^2M$ and {\sf
  Coker }$a$ (resp. {\sf Coker }$u$). From the
  commutative diagram:
$$\xymatrix{& \text{{\sf Coker }}a\ar[r]\ar[dd]^{\gamma} & X^2M\ar[r]\ar[dd] &
  SM\ar[rr]\ar@{=}[dd]\ar[rd]^{\Phi_M} & &\Sigma\text{{\sf Coker }}a\ar[dd]^{\Sigma\gamma}\\
  \Sigma^2M\ar[ur]^{\alpha}\ar[dr]^{\beta} & & & &
  \Sigma^3M\ar[ur]^{\Sigma\alpha} & \\ & \text{{\sf Coker }}u\ar[r] & R\ar[r] &
  SM\ar[rr]\ar[dr]^{\theta} & & \Sigma\text{{\sf Coker }}u \\ & & & &
  \Sigma^3M\ar[ur]^{\Sigma\beta}& }$$
we deduce the equality 
$\theta=(\Sigma\beta)^{-1}\gamma\Sigma\alpha\Phi_M=\Phi_M$ in
  $\underline{\modd}\cp$. The constructed $S$-complex is a $\Phi$-$S$-complex.
\\
 
\textbf{TR2}:
Let $\xymatrix{X:P\ar[r]^u & Q\ar[r]^v &R\ar[r]^w & SP}$ be a
$\Phi$-$S$-complex. It is homotopy-equivalent to a standard triangle
$T_M$. Thus the $S$-complex $$\xymatrix{X':Q\ar[r]^{-v} & R\ar[r]^{-w} &SP\ar[r]^{-Su} & SQ}$$
is homotopy-equivalent to $T_M[1]$. Since $\underline{T}$ is a
triangle functor, the objects $T_{\Sigma M}$ and $T_M[1]$ are
isomorphic in the stable category \underline{$S$-\comp}, \emph{i.e.}
they are homotopy-equivalent. Thus, by lemma \ref{lem0}, $T_M[1]$ is a
$\Phi$-$S$-complex and then so is $X'$.
\\

\textbf{TR3}:
This axiom is a direct consequence of lemma \ref{lem1}.
\\

\textbf{TR4}:
Let $X$ and $X'$ be two $\Phi$-$S$-complexes and suppose we have a
commutative diagram:
$$\xymatrix@1{X:P\ar[r]^(.6)u\ar@<2ex>[d]^{f^0} & Q\ar[r]^v\ar[d]^{f^1} &
  R\ar[r]^w & SP\ar[d]^{Sf^0}\\
X':P'\ar[r]^(.6){u'} & Q'\ar[r]^{v'} & R'\ar[r]^{w'} & SP'.}$$
Let $M$ (resp. $M'$) be the kernel of $u$ (resp. $u'$), and
  $g:M\rightarrow M'$ the
  induced morphism. The morphism $Tg:T_M\rightarrow T_{M'}$ induces a
  $S$-complex morphism $\tilde{g}=(g^0,g^1,g^2)$ between $X$ and
  $X'$. 

We are going to show that we can find a morphism $f^2:R\rightarrow R'$
such that $(f^0,f^1,f^2)$ can be extended in an $S$-complex morphism
that is homotopic to $\tilde{g}$. As $(g^0,g^1)$ and $(f^0,f^1)$
induce the same morphism $g$ in the kernels, we have some morphisms
$h^1:Q\rightarrow P'$ and $h^2:R\rightarrow Q'$ such that
$f^0-g^0=h^1u$ and $f^1-g^1=u'h^1+h^2v.$ We put $f^2=g^2+v'h^2$. We
have the following equalities:
$$\begin{array}{rclcrcl}
f^2v & = & g^2v+v'h^2v & \quad \textrm{and}\quad & w'f^2 & = & w'g^2
\\
 & = & v'(g^1+h^2v) & & & = & (Sg^0)w\\
 & = & v'(f^1-u'h^1) & & & = & (Sf^0-Sh^1Su)w\\
& = & v'f^1 & & & = & (Sf^0)w
\end{array}$$  
Thus $(f^0,f^1,f^2)$ can be extended to an $S$-periodic morphism
$\tilde{f}$ which is $S$-homotopic to $\tilde{g}$.
Their respective cones $C(\tilde{f})$ and $C(\tilde{g})$ are
isomorphic as $S$-complexes. Moreover, since $\tilde{g}$ is a
composition of $Tg:T_M\rightarrow T_{M'}$ with homotopy-equivalences,
the cones $C(\tilde{g})$ and $C(Tg)$ are homotopy-equivalent. 

In $\underline{\modd}\cp$, we have a triangle
$$\xymatrix{M\ar[r]^g &M'\ar[r] & C(g)\ar[r] & \Sigma M.}$$
Since $\underline{T}$ is a triangle functor, the sequence 
 $$\xymatrix{T_M\ar[r]^{T_g} &T_{M'}\ar[r] & T_{C(g)}\ar[r] &
 T_{\Sigma M}}$$
is a triangle in \underline{$S$-{\sf comp }}. But we know that 
 $$\xymatrix{T_M\ar[r]^{Tg} &T_{M'}\ar[r] & C(Tg)\ar[r] &
 T_M[1]}$$
is a triangle in $\underline{S\textrm{-\comp}}$.
Thus the objects $C(Tg)$ and $T_{C(g)}$ are isomorphic in
 $\underline{S\textrm{-\comp}}$, \emph{i.e.}
 homotopy-equivalent. Thus, the cone $C(\tilde{f})$ of $\tilde{f}$ is
 a $\Phi$-$S$-complex by lemma \ref{lem0}.

\section{Application to the deformed preprojective algebras}\label{defpreprojalg}

In this section, we apply the theorem \ref{thmprincipal} to show that
the category of finite dimensional projective modules over a deformed
preprojective algebra of generalized Dynkin type (see \cite{Bia}) is triangulated.
This will give us some examples of non standard triangulated categories
with finitely many indecomposables.
\subsection{Preprojective algebra of generalized Dynkin type}
Recall the notations of \cite{Bia}.
Let $\Delta$ be a generalized Dynkin graph of type $\mathbb{A}_n$,
$\mathbb{D}_n$ ($n\geq 4$), $\mathbb{E}_n$ ($n=6,7,8$), or
$\mathbb{L}_n$. Let $Q_\Delta$ be the following associated quiver:

$$\xymatrix{\Delta=\mathbb{A}_n \ (n\geq 1): & 0\ar@<.5ex>[r]^{a_0} &
  1\ar@<.5ex>[r]^{a_1}\ar@<.5ex>[l]^{\overline{a}_0} &
  2\ar@<.5ex>[l]^{\overline{a}_1} \ar@{.}[r] &
  n-2\ar@<.5ex>[r]^{a_{n-2}} &n-1\ar@<.5ex>[l]^{\overline{a}_{n-2}}}$$

$$\xymatrix{\Delta=\mathbb{D}_n \ (n\geq 4): & 0\ar@<.5ex>[dr]^{a_0} & &&&\\
   && 2\ar@<.5ex>[r]^{a_2}\ar@<.5ex>[ul]^{\overline{a}_0}
  \ar@<.5ex>[dl]^{\overline{a}_1} &
  3 \ar@<.5ex>[l]^{\overline{a}_2} \ar@{.}[r] &
  n-2\ar@<.5ex>[r]^{a_{n-2}} &n-1\ar@<.5ex>[l]^{\overline{a}_{n-2}} \\
  & 1\ar@<.5ex>[ur]^{a_1}  &&&& }$$

$$\xymatrix{\Delta=\mathbb{E}_n \ (n=6,7,8): & & &
  0\ar@<.5ex>[d]^{a_0} & & &\\
  & 1\ar@<.5ex>[r]^{a_1} &
  2\ar@<.5ex>[l]^{\overline{a}_1}\ar@<.5ex>[r]^{a_2} &
  3\ar@<.5ex>[l]^{\overline{a}_2}\ar@<.5ex>[r]^{a_3}\ar@<.5ex>[u]^{\overline{a}_0} &
  4\ar@<.5ex>[l]^{\overline{a}_3} \ar@{.}[r] &
  n-2\ar@<.5ex>[r]^{a_{n-2}} &n-1\ar@<.5ex>[l]^{\overline{a}_{n-2}}}$$

$$\xymatrix{\Delta=\mathbb{L}_n \ (n\geq 1): & & 0\ar@(ul,dl)_{\epsilon=\overline{\epsilon}}\ar@<.5ex>[r]^{a_0} &
  1\ar@<.5ex>[r]^{a_1}\ar@<.5ex>[l]^{\overline{a}_0} &
  2\ar@<.5ex>[l]^{\overline{a}_1} \ar@{.}[r] &
  n-2\ar@<.5ex>[r]^{a_{n-2}}
  &n-1\ar@<.5ex>[l]^{\overline{a}_{n-2}}}.$$
 
The \emph{preprojective algebra} $P(\Delta)$ associated to the graph
$\Delta$ is the quotient of the path algebra $kQ_\Delta$ by the relations:
$$\sum_{sa=i}a\overline{a},\quad\textrm{for each vertex $i$ of
  $Q_\Delta$.}$$
The following proposition is classical \cite[prop 2.1]{Bia}.
\begin{prop}
The preprojective algebra $P(\Delta)$ is finite dimensional and
selfinjective. Its Nakayama permutation $\nu$ is the identity for $\Delta=\mathbb{A}_1$, $\mathbb{D}_{2n}$,
$\mathbb{E}_7$, $\mathbb{E}_8$ and $\mathbb{L}_n$, and is of order $2$
in all other cases.
\end{prop}

\subsection{Deformed preprojective algebras of generalized Dynkin
  type}
Let us recall the definition of deformed preprojective algebra
introduced by \cite{Bia}. Let $\Delta$ be a graph of generalized
Dynkin type. We define an associated algebra $R(\Delta)$ as follows:
$$\begin{array}{rcl} R(\mathbb{A}_n)&=&k;\\ 
R(\mathbb{D}_n)&=&k\langle x,y\rangle /(x^2,y^2,(x+y)^{n-2});\\ 
R(\mathbb{E}_n)&=&k\langle x,y\rangle /(x^2,y^3,(x+y)^{n-3});\\
R(\mathbb{L}_n)&=&k[x]/(x^{2n}).
\end{array}$$
Further, we fix an exceptional vertex in each graph as follows (with
the notations of the previous section):
$$\begin{array}{rcl} 0&  \textrm{for} & \Delta=\mathbb{A}_n \textrm{
  or }
  \mathbb{L}_n,\\
2 & \textrm{for} & \Delta=\mathbb{D}_n,\\
3 & \textrm{for} & \Delta=\mathbb{E}_n.
\end{array}$$
Let $f$ be an element of the square $rad^2R(\Delta)$ of the radical of
$R(\Delta)$. The \emph{ deformed preprojective algebra} $P^f(\Delta)$
is the quotient of the path algebra $kQ_\Delta$ by the relations:  
$$\sum_{sa=i}a\overline{a},\quad \textrm{for each non exceptional vertex $i$ of $Q$,}$$
and $$\begin{array}{lcl} a_0\overline{a}_0 & \textrm{for}&\Delta=\mathbb{A}_n;\\
\overline{a}_0a_0+\overline{a}_1a_1+a_2\overline{a}_2+f(
  \overline{a}_0a_0,\overline{a}_1a_1),\textrm{ and } 
  (\overline{a}_0a_0+\overline{a}_1a_1)^{n-2} & \textrm{for}&\Delta=\mathbb{D}_n;\\
\overline{a}_0a_0+\overline{a}_2a_2+a_3\overline{a}_3+f(\overline{a}_0a_0,\overline{a}_2a_2),\textrm{
  and }
  (\overline{a}_0a_0+\overline{a}_2a_2)^{n-3} & \textrm{for}&\Delta=\mathbb{E}_n;\\
\epsilon^2+a_0\overline{a}_0 +\epsilon f(\epsilon), \textrm{ and }
  \epsilon^{2n} &
  \textrm{for}&\Delta=\mathbb{L}_n.\end{array}$$
Note that if $f$ is zero, we get the preprojective algebra $P(\Delta)$.

\subsection{Corollaries of ~\cite{Bia}}

The following proposition \cite[prop 3.4]{Bia} shows that the category
$\proj P^f(\Delta)$ of finite-dimensional projective modules over a
deformed preprojective algebra satisfies the hypothesis of theorem \ref{thmprincipal}.

\begin{prop}\label{BES}  
Let $A=P^f(\Delta)$ be a deformed preprojective algebra. Then there
exists an exact sequence of $A$-$A$-bimodules
$$\xymatrix{0\ar[r] & _1A_{\Phi^{-1}}\ar[r] & P_2\ar[r] & P_1\ar[r] &
  P_0\ar[r] & A\ar[r] & 0,}$$ where $\Phi$ is an automorphism of $A$
  and where the $P_i$'s are projective as bimodules. Moreover, for
  each idempotent $e_i$ of $A$, we have $\Phi(e_i)=e_{\nu(i)}$. 
\end{prop}   
So we can easily deduce the
corollary: 
\begin{cor}\label{cor1}
Let $P^f(\Delta)$ be a deformed preprojective algebra of generalized
Dynkin type. Then the category $\proj P^f(\Delta)$ of finite
dimensional projective modules is triangulated. The suspension is the
Nakayama functor. 
\end{cor}

Indeed, if $P_i=e_i A$ is a projective indecomposable, then
$P_i\ten_A A_\Phi$ is equal to $\Phi (e_i)A=e_{\nu (i)}A$ thus to $\nu (P_i)$.

Now we are able to answer to the question of the previous part and
find a triangulated category with finitely many indecomposables which
is not standard. The proof of the following theorem comes essentially from the theorem
\cite[thm 1.3]{Bia}.

\begin{thm}
Let $k$ be an algebraically closed field of characteristic $2$. Then
there exist $k$-linear triangulated categories with finitely many
indecomposables which are not standard.
\end{thm}

\begin{proof}
By theorem \cite[thm 1.3]{Bia}, we know that there exist basic deformed
preprojective algebras of generalized Dynkin type $P^f(\Delta)$ which
are not isomorphic to $P(\Delta)$. Thus the categories $\proj P^f(\Delta)$
and $\proj P(\Delta)$ can not be equivalent. But both are triangulated
by corollary \ref{cor1} and have the same AR-quiver
$\mathbb{Z}\Delta/\tau=Q_\Delta$.
\end{proof}

Conversely, we have the following theorem:

\begin{thm}
Let $\T$ be a finite 1-Calabi-Yau triangulated category. Then $\T$ is
equivalent to $\proj \Lambda$ as $k$-category, where
$\Lambda$ is a deformed preprojective algebra of generalized Dynkin type.
\end{thm}

\begin{proof}
Let $M_1, \ldots, M_n$ be representatives of the isoclasses of
indecomposable objects of $\T$. The $k$-algebra
$\Lambda=\End(\bigoplus_{i=1}^n M_i)$ is basic, finite-dimensional and selfinjective since $\T$
has a Serre duality. It is easy to see that $\T$ and $\proj \Lambda$
are equivalent as $k$-categories.

Let $\modd \Lambda$ be the category of finitely presented
$\Lambda$-modules. It is a Frobenius category. Denote by $\Sigma$ the
suspension functor of the triangulated category
$\underline{\modd}\Lambda$. The category $\T$ is $1$-Calabi-Yau, that
is to say that the suspension functor $S$ of the triangulated category
$\T$ and the Serre functor $\nu$ are isomorphic. But in
$\underline{\modd} \Lambda$, the functors $S$ and $\Sigma^3$ are
isomorphic. Thus, for each non projective simple $\Lambda$-module $M$
we have an isomorphism $\Sigma^3M\simeq \nu M$. By 
\cite[thm 1.2]{Bia}, we get immediately the result. 
\end{proof}

\bibliographystyle{amsplain}
\bibliography{biblio.bib}

\end{document}